\documentclass[12pt,notitlepage]{article}
\usepackage{latexsym,amsfonts,amssymb,amsmath,amsthm} 
\newcommand{\mz}{\ensuremath{\mathbb Z}}
\newcommand{\mr}{\ensuremath{\mathbb R}}

\newcommand{\mq}{\ensuremath{\mathbb Q}}
\newcommand{\mc}{\ensuremath{\mathbb C}}

\newcommand{\eop}{\ensuremath{\hfill \blacksquare}}

\newcommand{\mymod}{\ensuremath{\negthickspace \negmedspace \pmod}}
\newcommand{\shortmod}{\ensuremath{\negthickspace \negthickspace \negthickspace \pmod}}

\newcommand{\onehalf}{\ensuremath{ \frac{1}{2}}}
\newcommand{\notdiv}{\ensuremath{\not \; |}}
\newcommand{\notdivtext}{\ensuremath{\not |}}
\newcommand{\intR}{\int_{-\infty}^{\infty}}

\theoremstyle{plain}		
	\newtheorem{mytheo}{Theorem}[section]
	
	\newtheorem{myprop}[mytheo]{Proposition}
	\newtheorem{mycoro}[mytheo]{Corollary}
     \newtheorem{mylemma}[mytheo]{Lemma}
	\newtheorem{mydefi}[mytheo]{Definition}

\theoremstyle{remark}

\begin{document}

\title{Lower-Order Terms of the 1-Level Density of Families of Elliptic Curves}
\author{Matthew P. Young \footnote{This research was partially conducted during the period the author was employed by the Clay Mathematics Institute as a Liftoff Fellow. This research was partially supported by an NSF Mathematical Sciences Post-Doctoral Fellowship.} \\
American Institute of Mathematics \\
360 Portage Ave. \\
Palo Alto, CA 94306-2244 \\
myoung@aimath.org
}

\maketitle

\abstract{The Katz-Sarnak philosophy predicts that statistics of zeros of families of L-functions are strikingly universal.  However, subtle arithmetical differences between families of the same symmetry type can be detected by calculating lower-order terms of the statistics of interest.
In this paper we calculate lower-order terms of the 1-level density of some families of elliptic curves.  
We show that there are essentially two different effects on the distribution of low-lying zeros.  First, low-lying zeros  are more numerous in families of elliptic curves $E$ with relatively large numbers of points $\mymod{p}$.
Second, and somewhat surprisingly, a family with a relatively large number of primes of bad reduction has relatively fewer low-lying zeros.

We also show that the lower order term can grow arbitrarily large by taking a biased family with a relatively large number of points $\mymod{p}$ for all small primes $p$.
}  

\section{Introduction}
The Katz-Sarnak philosophy predicts that statistics of zeros of families of L-functions are strikingly universal.  However, subtle arithmetical differences between families of the same symmetry type can be detected by calculating lower-order terms of the statistics of interest.
In this paper we calculate lower-order terms of the 1-level density of some families of elliptic curves.  
We show that there are essentially two different effects on the distribution of low-lying zeros.  First, low-lying zeros  are more numerous in families of elliptic curves $E$ with relatively large numbers of points $\mymod{p}$.
Second, and somewhat surprisingly, a family with a relatively large number of primes of bad reduction has relatively fewer low-lying zeros.

We also show that the lower order term can grow arbitrarily large by taking a biased family with a relatively large number of points $\mymod{p}$ for all small primes $p$.
\subsection{Preliminaries and Notation}
Let $E/\mq$ be an elliptic curve with associated L-function $L(s, E)$ and conductor $N$.
By work of Wiles and others (\cite{Wiles}, \cite{TaylorWiles}, \cite{BCDT}) there is a weight two primitive holomorphic cusp form $f$ on $\Gamma_0(N)$ such that $L(s, E) = L(s, f)$.  The completed L-function $\Lambda(s, E)$ is entire and, with our normalization, satisfies the functional equation ${\Lambda(s, E) = \pm \Lambda(1-s, E)}$.

We are interested in the distribution of zeros of $\Lambda(s, E)$, especially zeros at or near the central point.  The Birch and Swinnerton-Dyer conjecture provides motivation for the interest in critical zeros, and random matrix theory provides tools for studying zeros near the central point.  It is our hope that insight can be gained into ranks of elliptic curves by studying the finer details of the distribution of zeros at or near the central point.

Consider the following `density,'
\begin{equation*}
D(E;\phi) = \sum_{\rho} \phi \left( \gamma \frac{\log{X}}{2\pi} \right), 
\end{equation*}
where $\phi$ is an even Schwartz-class function whose Fourier transform $\widehat{\phi}$ has compact support, $\rho = \onehalf + i \gamma$ runs over the nontrivial zeros of $L(s, E)$, and $X$ is a scaling parameter at our disposal.  Since $\phi$ has rapid decay, $D(E;\phi)$ is essentially a measure of zeros of distance $\ll (\log{X})^{-1}$ from the central point. In our notation we suppress the dependence of $D(E;\phi)$ on $X$.

\begin{mydefi}
Let $\mathcal{F}$ be a family of elliptic curves.
The {\em 1-level density} $\mathcal{D}_X(\mathcal{F})$ is by definition given by
\begin{equation*}
\mathcal{D}_X(\mathcal{F}) = \frac{1}{W_X(\mathcal{F})} \sum_{E \in \mathcal{F}}D(E;\phi) w_X(E),
\end{equation*}
where
\begin{equation*}
W_X(\mathcal{F}) = \sum_{E \in \mathcal{F}}w_X(E)
\end{equation*}
is the size of the family $\mathcal{F}$ and $w_X$ is a weighting function (to be thought of as choosing curves $E$ with conductor $\approx X$).
\end{mydefi}  

Our goal in this paper is to study the finer details of the asymptotics (as $X \rightarrow \infty$) of the 1-level density for some families of elliptic curves.  Katz and Sarnak have made predictions on the asymptotic behavior of the 1-level density for general families of L-functions (\cite{KSM}, \cite{KS}).  The asymptotic behavior is governed solely by the symmetry type of the family.  However, the lower-order terms in the 1-level density heavily depend on the arithmetical nature of the family.

There is a large body of work on the 1-level density for a variety of families of L-functions (e.g. \cite{FI}, \cite{ILS}, \cite{KSM}, \cite{R}).  Additionally, S. J. Miller \cite{M} and the author \cite{Young} have investigated the main term of the 1-level density for various families of elliptic curves.  Miller concentrates upon one-parameter families and considers both the 1-level and 2-level densities.  The author focuses on large two-parameter families with test functions $\phi$ concentrated close to the origin.  Fouvry and Iwaniec \cite{FI} also obtained lower-order terms in their investigations.

Conrey and Snaith \cite{CS} have developed a method for predicting lower order terms for an arbitrary family of L-functions.  This method is based on the L-functions ratios conjectures of Conrey, Farmer, and Zirnbauer \cite{CFZ}.  Their work is powerful in its generality and conceptually quite illuminating, but of course it does not give unconditional results.  If one were to apply their method to one of the families considered in this paper, then one would be led to essentially the same (unconditional) calculations that we carry out in this paper.


\subsection{Statement of Results}
The families of elliptic curves under scrutiny in this paper are the family of curves with torsion group $\mz/2\mz \times \mz/2\mz$ and the family with torsion group $\mz/2\mz$.  We study these families (as opposed to the family of all elliptic curves) because with these families we can more accurately estimate the average of the logarithm of the conductors (see Lemmas \ref{lem:conductorestimate} and \ref{lem:conductorestimate2}).

The following definitions precisely state which families we are considering and what are the corresponding weighting functions $w_X$.

\begin{mydefi}
Let $q$ be a positive odd integer, and let $a_0$ and $b_0$ be integers such that $(q, a_0 b_0 (a_0 + 2 b_0)) = 1$.  The family $\mathcal{F}_1 = \mathcal{F}_1(a_0, b_0 ;q) = \{E_{a, b} \}$ is the family of elliptic curves given by the Weierstrass equations $E_{a, b}:y^2 = x(x-a)(x+2b)$, where $a$ and $b$ are odd coprime integers such that $a \equiv a_0 \pmod{q}$ and $b \equiv b_0 \pmod{q}$.  Further, let $N_{a, b}$ be the conductor of $E_{a, b}$ and let $\lambda_{a, b}(n)$ be the coefficient of $n^{-s}$ in the Dirichlet series expansion of $L(s, E_{a, b})$.  In the special case $q=1$ we set
\begin{equation*}
\mathcal{F}_1' = \mathcal{F}_1(1, 1;1).
\end{equation*}
\end{mydefi}
Any elliptic curve over $\mq$ with torsion group $\mz/2\mz \times \mz/2\mz$ can be represented by a Weierstrass equation of the form $y^2 = x(x-a)(x+b)$.  We shall refer to $\mathcal{F}_1'$ as {\em the} family of all elliptic curves with torsion $\mz/2\mz \times \mz/2\mz$, even though the coprimality restrictions in the definition of $\mathcal{F}_1'$ mean that our family is somewhat smaller.
The coprimality restrictions are imposed for technical reasons, namely, in order to accurately compute the conductor; see Lemma \ref{lem:conductor} for a formula for $N_{a, b}$.  

It is interesting to study the variation of the 1-level density of $\mathcal{F}_1$ as $a_0$ and $b_0$ vary $\pmod{q}$.  
The coprimality condition $(q, a_0 b_0 (a_0 + 2b_0)) = 1$ means $E_{a, b}$ has good reduction at every prime dividing $q$.  

The corresponding definition for the family with torsion $\mz/2\mz$ is given by
\begin{mydefi}
Let $q$ be a positive odd integer, and let $a_0$ and $b_0$ be integers such that $(q, b_0 (a_0^2 +  b_0)) = 1$.  The family $\mathcal{F}_2 = \mathcal{F}_2(a_0, b_0 ;q) = \{E_{a, b} \}$ is the family of elliptic curves given by the Weierstrass equations $E_{a, b}:y^2 = x(x^2 + 2ax -b)$, where $a$ and $b$ are coprime integers such that $a \equiv a_0 \pmod{q}$, $b \equiv b_0 \pmod{q}$, and $a \equiv b \equiv 1 \pmod{4}$.  Further, let $N_{a, b}$ be the conductor of $E_{a, b}$ and let $\lambda_{a, b}(n)$ be the coefficient of $n^{-s}$ in the Dirichlet series expansion of $L(s, E_{a, b})$.
In the special case $q=1$ we set
\begin{equation*}
\mathcal{F}_2' = \mathcal{F}_2(1, 1;1).
\end{equation*}
\end{mydefi}

For notational simplicity we use the common symbols $E_{a, b}$, $N_{a, b}$, and $\lambda_{a, b}(n)$ for both families $\mathcal{F}_1$ and $\mathcal{F}_2$.  There should be no confusion because we shall always make clear to which family we are referring.

We define $w_X$ for $\mathcal{F}_1$ with the following

\begin{mydefi}
Let $w \in C_0^{\infty} (\mr^+ \times \mr^+)\footnote{Here and throughout $\mr^+ = (0, \infty)$.}$, $\widehat{w}(0, 0) = 1$.  Set $A = B = X^{1/3}$.  Then for $E_{a, b} \in \mathcal{F}_1$ we set
\begin{equation*}
w_X(E_{a, b}) = w \left( \frac{a}{A}, \frac{b}{B} \right).
\end{equation*}
\end{mydefi}
Likewise, for $\mathcal{F}_2$ we have
\begin{mydefi}
Let $w \in C_0^{\infty} (\mr^+ \times \mr^+)$, $\widehat{w}(0, 0) = 1$.  Set $A = X^{1/4}, B = X^{1/2}$.  Then for $E_{a, b} \in \mathcal{F}_2$ we set
\begin{equation*}
w_X(E_{a, b}) = w \left( \frac{a}{A}, \frac{b}{B} \right).
\end{equation*}
\end{mydefi}

The conditions on the sizes of $A$ and $B$ are `correct', because in Section \ref{section:conductor} we show that for both families $\mathcal{F}_i$
\begin{equation*}
\sum_{E_{a, b} \in \mathcal{F}_i} \frac{\log{N_{a, b}}}{\log{X}} w_X(E_{a, b}) \sim W_X (\mathcal{F}_i) \; \; \; \text{as $X \rightarrow \infty$},
\end{equation*}
so that we can think of $w_X$ as picking out curves $E$ with conductor $N \approx X$.  The reason that these are the right choices for $A$ and $B$ is that, for $\mathcal{F}_1$, say, $N_{a, b} = 2^5 a^* b^* (a + 2b)^*$ (here $n^*$ is the product of primes dividing $n$).  If $w_X(E_{a, b}) \neq 0$, then $a \asymp A$ and $b \asymp B$, so for `most' values of $a$ and $b$, we should have $\log{N_{a, b}}$ close to $\log(AB(A + 2B)) \sim \log{X}$.  Similar reasoning holds for $\mathcal{F}_2$.  

Now we may state our main theorems
\begin{mytheo}
\label{thm:mainthmloworder}
Let $q$ be an odd positive integer. Then
\begin{align*}
\mathcal{D}_X(\mathcal{F}_1) & =    \frac{1}{2} \phi(0) + \widehat{\phi}(0) 
+ \frac{\widehat{\phi}(0)}{\log{X}} \left[e_1(a_0, b_0) + d_{1, 1}(q) + d_{3, 1}(q) +d_{4, 1}(q) + d_{5, 1}(q) + d_{6, 1}(q) \right]\\
  & \qquad + \frac{\widehat{\phi}(0)}{\log{X}} \left[ c_{1, 1} + c_{2, 1} + c_{3, 1} + c_{4, 1} + c_{5, 1} + c_{6, 1} \right] + O\left(\frac{1}{\log^3{X}}\right)
\end{align*}
as $X \rightarrow \infty$, provided supp $\widehat{\phi} \subset (-\frac{2}{3}, \frac{2}{3})$, and where the constants are given by \eqref{eq:d1}, \eqref{eq:c1}, \eqref{eq:c2}, \eqref{eq:d3}, \eqref{eq:c3}, \eqref{eq:d4}, \eqref{eq:c4}, \eqref{eq:c5}, \eqref{eq:c6}, \eqref{eq:d5}, \eqref{eq:d6}, and \eqref{eq:e1}.
\end{mytheo}
For convenience to the reader we gather the constants together here.  Let $\chi_{N_{a,b}}$ be the principal Dirichlet character $\mymod{N_{a, b}}$.  We have
\begin{align*}
e_1(a_0, b_0) & = -2 \sum_{p | q} \log{p} \left(1 - \frac1p \right) \left\{\left(1 - \frac{\lambda_{a_0, b_0}(p)}{p^{1/2}} + \frac{1}{p} \right)^{-1} -1 \right\}, \\
d_{1, 1}(q) & = 3 \sum_{p | q} \frac{\log{p}}{p^2 - 1}, \\
d_{3, 1}(q) & = 6 \sum_{p | q} \frac{ \log{p}}{(p^2 - 1)(p + 1)}, \\
d_{4, 1}(q) & = 6 \sum_{p | q} \frac{\log{p}}{p(p+1)}, \\
d_{5, 1}(q) & = 2 \sum_{p | q} \frac{\log{p}}{p(p + 1)} \mathop{\sum \sum}_{a, b \shortmod{p}} \left\{ \left(1 - \frac{\lambda_{a, b}(p)}{p^{1/2}} + \frac{\chi_{N_{a,b}}(p)}{p} \right)^{-1} - 1\right\}, \\
& \text{and} \\
d_{6, 1} (q)  & = - 6 \sum_{p | q} \frac{\log{p}}{p(p+1)^2} 
\end{align*}
Note that $\chi_{N_{a, b}}(p)$ is well-defined in the above expression for $d_{5, 1}(q)$, even though $N_{a, b}$ is not (the point is that the property of $p$ dividing the conductor $N_{a, b}$ of $y^2 = x(x-a)(x+2b)$ only depends on $a$ and $b \pmod{p}$).

As for the $c_{i, 1}$ we have
\begin{align*}
c_{1, 1} & =  \int_0^\infty \int_0^\infty \log(2^5 xy(x+2y)) w(x, y) dx dy - 3 \sum_{p \neq 2} \frac{\log{p}}{p^2 - 1}, \\
c_{2, 1} & = -2 \log{2\pi} - 2 \gamma, \\
c_{3, 1} & = -6 \sum_{p \neq 2} \frac{ \log{p}}{(p^2 - 1)(p + 1)}, \\
c_{4, 1} & =  2 \left(1 + \int_1^{\infty} \frac{R(t)}{t^2} dt  - 3 \sum_{p \neq 2} \frac{\log{p}}{p(p+1)} - \log{2} \right), \\
c_{5, 1} & =  -2   \sum_{p \neq 2} \frac{\log{p}}{p(p + 1)} \mathop{\sum \sum}_{a, b \shortmod{p}} \left\{ \left(1 - \frac{\lambda_{a, b}(p)}{p^{1/2}} + \frac{\chi_{N_{a,b}}(p)}{p} \right)^{-1} - 1\right\}, \\
& \text{and} \\
c_{6, 1} & =  6 \sum_{p \neq 2} \frac{\log{p}}{p(p + 1)^2},
\end{align*}
and where $R(t) = \sum_{p \leq t} \log{p}  - t$.  Actually, the formulas given above for $e_1, d_{5, 1}$, and $c_{5,1}$ are given by (\ref{eq:e1}'), (\ref{eq:d5}'), and (\ref{eq:c5}'), respectively.

It may be of interest that the calculation of the lower-order terms is unconditional, but that obtaining the result for support up to $2/3$ relies on the Riemann Hypothesis for Dirichlet L-functions (but not the elliptic curve L-functions under consideration!).  Eliminating the use of the GRH would force us to reduce the support.  We also have
\begin{mycoro}
\label{coro:loworder}
We have
\begin{eqnarray*}
\mathcal{D}_X(\mathcal{F}_1')  =   \frac{1}{2} \phi(0) + \widehat{\phi}(0) 
  + \frac{\widehat{\phi}(0)}{\log{X}} \left[c_{1, 1} + c_{2, 1} + c_{3, 1} + c_{4, 1} + c_{5, 1} + c_{6, 1} \right] + O\left(\frac{1}{\log^3{X}}\right)
\end{eqnarray*}
as $X \rightarrow \infty$, provided supp $\widehat{\phi} \subset (- \frac{2}{3}, \frac{2}{3})$, and where the constants are as in Theorem \ref{thm:mainthmloworder}.
\end{mycoro}

Similarly, for the family $\mathcal{F}_2$ we have
\begin{mytheo}
\label{thm:mainthmloworder2}
Let $q$ be an odd positive integer.  Then
\begin{align*}
\mathcal{D}_X(\mathcal{F}_2) & =    \frac{1}{2} \phi(0) + \widehat{\phi}(0) 
+ \frac{\widehat{\phi}(0)}{\log{X}} \left[e_2(a_0, b_0) + d_{1, 2}(q) + d_{2, 2}(q) + d_{3, 2}(q) + d_{4, 2}(q) + d_{5, 2}(q) \right]\\
& \qquad  + \frac{\widehat{\phi}(0)}{\log{X}} \left[ c_{1, 2} + c_{2, 2} + c_{3, 2} + c_{4, 2} + c_{5, 2} + c_{6, 2}  \right] + O\left(\frac{1}{\log^3{X}}\right)
\end{align*}
as $X \rightarrow \infty$, provided supp $\widehat{\phi} \subset (-\frac{1}{2}, \frac{1}{2})$, and where the constants are given by \eqref{eq:d1'}, \eqref{eq:c1'}, \eqref{eq:c2}, \eqref{eq:d3'}, \eqref{eq:c3'}, \eqref{eq:d4}, \eqref{eq:c4}, \eqref{eq:c5}, \eqref{eq:c6}, \eqref{eq:d5}, \eqref{eq:d6} and \eqref{eq:e1}. 
\end{mytheo}
Again we presently reproduce the constants.  We have
\begin{align*}
e_2(a_0, b_0) & = -2 \sum_{p | q} \log{p} \left(1 - \frac1p \right) \left\{ \left(1 - \frac{\lambda_{a_0, b_0}(p)}{p^{1/2}} + \frac{1}{p} \right)^{-1} -1 \right\}, \\
d_{1, 2}(q) & = 2 \sum_{p | q} \frac{\log{p}}{p^2 - 1}, \\
d_{3, 2}(q) & = 4 \sum_{p | q} \frac{ \log{p}}{(p^2 - 1)(p + 1)}, \\
d_{4, 2}(q) & = 4 \sum_{p | q} \frac{\log{p}}{p(p+1)}, \\
d_{5, 2}(q) & = 2 \sum_{p | q} \frac{\log{p}}{p(p + 1)} \mathop{\sum \sum}_{a, b \shortmod{p}} \left\{ \left(1 - \frac{\lambda_{a, b}(p)}{p^{1/2}} + \frac{\chi_{N_{a,b}}(p)}{p} \right)^{-1} - 1\right\}, \\
& \text{and} \\
d_{6, 2} (q)  & = - 4 \sum_{p | q} \frac{\log{p}}{p(p+1)^2} 
\end{align*}
For the $c_{i, 2}$ we have
\begin{align*}
c_{1, 2} & =  \int_0^\infty \int_0^\infty \log(2^6 y(x^2+y)) w(x, y) dx dy - 2 \sum_{p \neq 2} \frac{\log{p}}{p^2 - 1}, \\
c_{2, 2} & = -2 \log{2\pi} - 2 \gamma, \\
c_{3, 2} & = -4 \sum_{p \neq 2} \frac{ \log{p}}{(p^2 - 1)(p + 1)}, \\
c_{4, 2} & =  2 \left(1 + \int_1^{\infty} \frac{R(t)}{t^2} dt  - 2 \sum_{p \neq 2} \frac{\log{p}}{p(p+1)} - \log{2} \right), \\
c_{5, 2} & =  -2   \sum_{p \neq 2} \frac{\log{p}}{p(p + 1)} \mathop{\sum \sum}_{a, b \shortmod{p}} \left\{ \left(1 - \frac{\lambda_{a, b}(p)}{p^{1/2}} + \frac{\chi_{N_{a,b}}(p)}{p} \right)^{-1} - 1\right\}, \\
& \text{and} \\
c_{6, 2} & =  4 \sum_{p \neq 2} \frac{\log{p}}{p(p + 1)^2},
\end{align*}
and where $R(t) = \sum_{p \leq t} \log{p}  - t$.

Here the support range is unconditional.  It is likely that support up to $2/3$ could be reached using methods in \cite{Young}, conditional on the GRH.

Again, we have
\begin{mycoro}
\label{coro:loworder2}
The following holds
\begin{eqnarray*}
\mathcal{D}_X(\mathcal{F}_2')  =   \frac{1}{2} \phi(0) + \widehat{\phi}(0) 
+ \frac{\widehat{\phi}(0)}{\log{X}} 
   \left[ c_{1, 2} + c_{2, 2} + c_{3, 2} + c_{4, 2} + c_{5, 2} + c_{6, 2} \right] + O\left(\frac{1}{\log^3{X}}\right)
\end{eqnarray*}
as $X \rightarrow \infty$, provided supp $\widehat{\phi} \subset (- \frac{1}{2}, \frac{1}{2})$, and where the constants $c_{i, 2}$ are as in Theorem \ref{thm:mainthmloworder2}.
\end{mycoro}

It is of interest to see how the 1-level densities of the families $\mathcal{F}_1$ and $\mathcal{F}_2$ vary as $a_0$ and $b_0$ vary $\mymod{q}$.  Since $e_i$ is the only constant depending on $a_0$ and $b_0$, Theorems \ref{thm:mainthmloworder} and \ref{thm:mainthmloworder2} show that the conditions $a \equiv a_0 \pmod{q}$ and $b \equiv b_0 \pmod{q}$ simply fix $\lambda_{a, b}(p) = \lambda_{a_0, b_0}(p)$ for $p | q$ (this can be seen as a sort of independence of primes since restricting $a$ and $b \pmod{q}$ does not affect the average behavior of $\lambda_{a, b}(p)$ for $(p, q) = 1$).  We clearly see how relatively large $-\lambda_{a_0, b_0}(p)$ corresponds to (slightly) more zeros at or near the central point.  It is well-known that an elliptic curve with positive rank should have larger $-\lambda(p)$'s than an elliptic curve with zero rank; here we see a direct relationship between low-lying zeros and larger $-\lambda(p)$'s.  Conrey {\it et al} \cite{CKRS} have noticed similar phenomena in the distribution of rank frequencies of quadratic twists of a fixed L-function.
Namely, they observed a similar kind of independence of primes and also a slight increase in rank frequency when restricting to arithmetic progressions with a larger $-\lambda(p)$.

In Theorems \ref{thm:mainthmloworder} and \ref{thm:mainthmloworder2} we can also let $q$ grow (slowly) with respect to $X$.  In our proofs of the two theorems we exhibit the dependence on $q$ of the implied constants in the remainder term.  As long as $q \ll X^{\varepsilon}$ for $\varepsilon$ small enough with respect to the support of $\widehat{\phi}$ we retain the desired asymptotic.  The payoff is that we may choose a sequence of $q$'s with biased $\lambda_{a_0, b_0}(p)$ for all $p | q$ so that $e_i(a_0, b_0)$ can grow with $X$.  Precisely we let $q_n$ be the product of all odd primes $p$ less than or equal to $n$.  For each $p$ select integers $a_p$ and $b_p$ such that $E_{a_p, b_p}$ has good reduction at $p$.  By the Chinese Remainder Theorem take $a_{0, n}$ and $b_{0,n}$ such that $a_{0,n} \equiv a_p \pmod{p}$ and $b_{0,n} \equiv b_p \pmod{p}$ for all $p \leq n$.  Notice that since $E_{a_p, b_p}$ has good reduction at all odd primes $p \leq n$, the conductor of each curve in the family is almost prime (that is, the conductor has no small prime factors except for the fixed power of $2$).
Then we may state
\begin{mytheo}
\label{thm:qvaries}
Take $\varepsilon$ small and let $n = n(X)$ be an integer such that
\begin{equation*}
n \sim \varepsilon \log{X}
\end{equation*}
as $X \rightarrow \infty$.
Then for $i=1, 2$ we have
\begin{equation*}
\mathcal{D}_X(\mathcal{F}_i(a_{0, n}, b_{0,n}; q_{n})) = \frac{1}{2} \phi(0) + \widehat{\phi}(0) + \frac{\widehat{\phi}(0)}{\log{X}} e_i(a_{0,n}, b_{0,n}) + C_i \frac{1}{\log{X}} + O\left( \frac{1}{\log^2{X}} \right)
\end{equation*}
as $X \rightarrow \infty$, provided $\widehat{\phi}$ has support as in Theorems \ref{thm:mainthmloworder} and \ref{thm:mainthmloworder2}.  Here $C_i$ is a fixed constant depending only on $w$, $\phi$, and the family $\mathcal{F}_i'$, and $\varepsilon$ is chosen small enough with respect to the support of $\widehat{\phi}$.
\end{mytheo}
For convenience we recall
\begin{align*}
e_i(a_{0,n}, b_{0,n}) & = -2 \sum_{p \leq n} \log{p} \left(1 - \frac{1}{p} \right) \left\{ \left(1 - \frac{\lambda_{a_p, b_p}(p)}{p^{1/2}} + \frac{1}{p} \right)^{-1} -1 \right\} \\
& = -2 \sum_{p \leq n} \log{p} \left( \frac{\lambda_{a_p, b_p}(p)}{p^{1/2}} - \frac{1}{p} + \frac{\lambda^2_{a_p, b_p}(p)}{p} \right) + O(1) \\
& = 2 \log{n} -2 \sum_{p \leq n} \log{p} \frac{\lambda_{a_p, b_p}(p)}{p^{1/2}} \left( 1 + \frac{\lambda_{a_p, b_p}(p)}{p^{1/2}} \right) +o(\log{n}).
\end{align*}

By taking biased $\lambda_{a_p, b_p}(p)$'s we can show
\begin{align*}
e_i(a_{0,n}, b_{0,n}) & \gg (\log{X})^{1/2}
\end{align*}
and for a different choice of $a_p$, $b_p$'s we have
\begin{align*}
e_i(a_{0,n}, b_{0,n}) & \ll -(\log{X})^{1/2}.
\end{align*}
These estimates translate to
\begin{equation*}
\mathcal{D}_X(\mathcal{F}_i(a_{0, n}, b_{0,n}; q_{n})) \geq \frac{1}{2} \phi(0) + \widehat{\phi}(0) + C_i \frac{1}{(\log{X})^{1/2}}
\end{equation*}
and 
\begin{equation*}
\mathcal{D}_X(\mathcal{F}_i(a_{0, n}, b_{0,n}; q_{n})) \leq \frac{1}{2} \phi(0) + \widehat{\phi}(0) - C_i \frac{1}{(\log{X})^{1/2}},
\end{equation*}
respectively, where $C_i$ is some positive constant.  These estimates show there is a lot of variation in the low-lying zero density.

We now show why there exist such large values of $e_i$.  It will follow from the fact that
\begin{equation*}
\text{max}_{a, b \shortmod{p}} |\lambda_{a,b}(p)| \geq 1 + O(p^{-1})
\end{equation*}
and the elementary fact that for any $\lambda_{a,b}(p)$ there exists $a'$ and $b'$ such that $\lambda_{a', b'}(p) = - \lambda_{a, b}(p)$.  The above estimate on the maximal size of $\lambda(p)$ is immediately deduced from
\begin{equation*}
\sum_{a \shortmod{p}} \sum_{b \shortmod{p}} \lambda^2_{a,b}(p) = p^2 + L_i(p),
\end{equation*}
where $L_i$ is a linear polynomial depending on the family only.  This elementary computation is performed in Section \ref{section:complete}

It is important to realize that we cannot take biased $\lambda(p)$'s for all $p \leq n$ and $n$ too large because it would necessarily contradict the Riemann hypothesis for an elliptic curve.  Nevertheless, there should be quite a bit of freedom in choosing $\lambda(p)$'s for $p$ relatively small.  

Notice also that we are heavily using the fact that our elliptic curves are given by Weierstrass equations in order to piece together our family from local data.  It would be of interest to make an analogous construction for general holomorphic cusp forms of, say, weight 2.  Loosely speaking, suppose we are given a finite set of primes $\{p_i \}$ and modular L-functions $\{L(s, f_i)\}$ of weight 2 and relatively small level.  Then the question is, can we produce a `natural' family of modular L-functions $\{L(s, f)\}$ (of weight 2 and ostensibly larger level) such that $\lambda_f(p_i) = \lambda_{f_i}(p_i)$ for all $f$ and $p_i$?

\subsection{Analysis of the Constants and a Comparison of the Two Families}
It is of interest to compare the lower-order terms of $\mathcal{F}_1'$ and $\mathcal{F}_2'$ (comparing general $\mathcal{F}_1(a_0, b_0, q)$ and $\mathcal{F}_2(a_0', b_0', q')$ is not straightforward since it is not clear how the congruence conditions on $\mathcal{F}_1$ should correspond to the congruence conditions on $\mathcal{F}_2$).  Recall that the family $\mathcal{F}_1'$ corresponds to a family with torsion $\mz/2\mz \times \mz/2\mz$; the family $\mathcal{F}_2'$ corresponds to a family with torsion $\mz/2\mz$.

In order to compare the lower-order terms of the two families we must scrutinize the constants $c_{i,j}$.  It is important to know what quantities are averaged to give these numbers.  The interested reader should see Proposition \ref{prop:explicitformula} for the form in which we have written the explicit formula.  There are five distinct quantities over which we shall average; the constants $c_{1,i}$ through $c_{4, i}$ correspond to the first four such quantities, respectively.  Averaging over the quantity on the last line of Proposition \ref{prop:explicitformula} gives both $c_{5,i}$ and $c_{6, i}$; we remove the restriction $p \notdivtext N$ by summing over $p | N$ separately.  Summing over $p|N$ gives $c_{6, i}$, and summing over all $p$ gives $c_{5, i}$.

The first thing to notice after examining the various constants is that $c_{1,i}$ is the only constant associated to $\mathcal{F}_i$ that depends on the weighting function $w_X$.  We must properly choose our weighting functions in order to compare the two families in a natural way.  We contend that the correct choices of $w_X$ should equalize $c_{1,1}$ and $c_{1,2}$.  The constants $c_{1,i}$ give the lower-order approximation to the average of the logarithm of the conductors for the two families.  The total number of zeros in the critical strip of height up to $T$ of an elliptic curve L-function with conductor $N$ is
\begin{equation*}
\frac{T}{\pi} \log{\frac{NT^2}{(2\pi e)^2}} + O(\log{NT}).
\end{equation*}
See \cite{IK}, Theorem 5.8, for instance.  Therefore, equalizing $c_{1,1}$ and $c_{1,2}$ is equivalent to making the proper normalization of the high zeros of the two families.

Obviously $c_{2,1} = c_{2,2}$ because all elliptic curve L-functions have the same gamma factors in their functional equation.

Computing the numerical values of the remaining $c_{i,j}$'s is enlightening.  We have\footnote{I believe all constants in this section are calculated correctly to the displayed number of digits, but I have not rigorously bounded the tails of the infinite series.}
\begin{align*}
c_{3,1} & = -0.3309763 \ldots \\
c_{3,2} & = -0.2206508 \ldots.
\end{align*}
The term $c_{3,i}$ is a measure of the effect of the primes dividing the conductor.  For $\nu$ odd the quantity $\lambda_{a, b}^\nu(p)$ is oscillatory and does not contribute even a lower order term upon averaging over the family.  On the other hand, $\lambda_{a, b}^{2\nu}(p) = p^{-\nu}$ for primes dividing the conductor, so the effect we are measuring is {\em how many} primes divide the conductor.  The difference is that $\mathcal{F}_1$ has three essentially independent large factors and $\mathcal{F}_2$ has two essentially independent large factors (precisely, we show in Section \ref{section:invariants} that for $\mathcal{F}_1$ we have $N_{a, b} = 2^5a^* b^*(a + 2b)^*$ and for $\mathcal{F}_2$ we have $N_{a, b} = 2^6 b^* (a^2 + b)^*$, $n^*$ here being the product of primes dividing $n$).

Next in this list we have
\begin{align*}
c_{4,1} & = -3.6656429 \ldots \\
c_{4,2} & = -3.127581 \ldots.
\end{align*}
Again, the difference between the two families is the number of primes dividing the conductor.

The most arithmetically interesting part arises from averaging $\lambda(p^\nu)$.  Here we compute
\begin{align*}
c_{5,1} & = -0.169117 \ldots \\
c_{5,2} & = -0.000614 \ldots.
\end{align*}
The reason $c_{5,2}$ is significantly smaller than $c_{5,1}$ is that for $\mathcal{F}_2$ we have
\begin{equation*}
\sum_{a \shortmod{p}} \sum_{b \shortmod{p}} \lambda_{a, b}(p^4) = 0
\end{equation*}
for all $p$, whereas the corresponding sum for $\mathcal{F}_1$ is nonzero in general (numerical computations suggest that it is zero for $p =2$ only).  We also have
\begin{align*}
c_{6,1} & = 0.24266089 \ldots \\
c_{6,2} & = 0.16177739 \ldots.
\end{align*}
Here again we have $c_{6,i}$ proportional to the number of independent divisors of the conductor of the family $\mathcal{F}_i$.

By adding the relevant constants we see
\begin{align*}
c_{3,1} + c_{4, 1} + c_{5,1} + c_{6,1} & = -3.77087 \ldots \\
c_{3,2} + c_{4, 2} + c_{5,2} + c_{6,2} & = -3.18707 \ldots,
\end{align*}
so that $\mathcal{F}_2$ has more low-lying zeros than $\mathcal{F}_1$.  It is somewhat surprising that the lower-order terms for these two families are controlled more by the number of prime divisors of the conductors than by any subtle distributional properties of the $\lambda_{a, b}(p^\nu)$.  Of course, for any $p$ prime and $s \geq 1/2$ we have
\begin{equation*}
\left(1 - \frac{\lambda(p)}{p^s} + \frac{1}{p^{2s}}\right)^{-1} < \left(1 - \frac{\lambda(p)}{p^s}\right)^{-1},
\end{equation*}
so, all things otherwise being equal, there should be a correlation between many primes dividing the conductor $N$ of an elliptic curve $E$ and slightly larger values of $L(1/2, E)$ (and perhaps slightly fewer low-lying zeros of $L(s, E)$).

\subsection{Overview of the Proofs and Organization of the Paper}
The proofs of Theorems \ref{thm:mainthmloworder} and \ref{thm:mainthmloworder2} follow from a careful evaluation using the explicit formula.  The two proofs are exceedingly similar so we have endeavored to eliminate repetition as much as possible.  In fact, we prove the two theorems in parallel.
As a rule, we do the necessary computations for $\mathcal{F}_1$ and $\mathcal{F}_2$ in separate subsections of the same section.  Often we first do a calculation more general than what we need and simply apply it to both families of interest.  It is usually no extra work to do the calculations in greater generality and in fact it reduces repetition.

The result of Corollary \ref{coro:loworder} (resp. Corollary \ref{coro:loworder2}) will follow from the proof of Theorem \ref{thm:mainthmloworder} (resp. Theorem \ref{thm:mainthmloworder2}); the calculations will go through with $q=1$.  The constants $e_i$, $d_{1, i}$, and $d_{2, i}$ will be zero in this case.  Of course the restrictions $a \equiv a_0 \pmod{q}$ etc. are then always satisfied for $q=1$.

The proof of Theorem \ref{thm:qvaries} follows from the same lemmas that prove Theorems \ref{thm:mainthmloworder} and \ref{thm:mainthmloworder2}.  The uniformity proven with respect to $q$ allows $q \ll X^{\varepsilon}$.  Taking $q$ to be the product of the odd primes less than or equal to $n$ means we must restrict $n$ by $n\ll \varepsilon \log{X}$.  One minor issue that arises is describing the behavior of the constants $d_{l, i}(q)$ as $q \rightarrow \infty$.  It is easily seen these converge to the sum over primes in the corresponding $c_{l, i}$, with error term of size $O(n^{-1}) = O((\log{X})^{-1})$.

We state the explicit formula for an L-function attached to an elliptic curve in Section \ref{section:explicitformula} (see Proposition \ref{prop:explicitformula}).  The overall goal, then, is to average each of the terms in the explicit formula.  We average each term in a separate section.  

In Section \ref{section:invariants} we calculate the conductor of the curves $E_{a,b}$ and the Dirichlet coefficients of $L(s, E_{a, b})$.

The conditions defining the families $\mathcal{F}_1$ and $\mathcal{F}_2$ are slightly thorny.  In Section \ref{section:sizes} we compute the sizes of the families, that is, we calculate the asymptotics of $W_X(\mathcal{F}_i)$.

In the remaining sections we average the various terms of the explicit formula.  The difficult parts are averaging the logarithm of the conductors (Section \ref{section:conductor}) and computing the sums in Section \ref{section:variation}.

\subsection{Acknowledgements}
I would like to thank Henryk Iwaniec and Brian Conrey for valuable advice.  Part of the material in this paper originally appeared in my PhD thesis under the advisement of Professor Iwaniec, who suggested this line of research.

\section{The Explicit Formula}
\label{section:explicitformula}
We need a precise formulation of the explicit formula.
Let $E$ be an elliptic curve with L-function 
\begin{equation*}
L(s, E) = \sum_{n=1}^{\infty} \frac{\lambda(n)}{n^s} = \prod_{p} \left(1 - \frac{\lambda(p)}{p^{s}} + \frac{\chi_0(p)}{p^{2s}} \right)^{-1},
\end{equation*}
where $\chi_0$ is the principal Dirichlet character modulo the conductor $N$ of $E$.  We normalize $L(s, E)$ to have central point $s = 1/2$.
The completed L-function 
\begin{equation*}
\Lambda(s, E) = \left( \frac{\sqrt{N}}{2\pi} \right)^{s + \onehalf} \Gamma(s + \onehalf) L(s, E)
\end{equation*}
is entire and satisfies the functional equation $\Lambda(s, E) = w \Lambda(1-s, E)$ where $w = \pm 1$ is the root number of $E$.

Now we simply apply Theorem 5.12 of \cite{IK} to $L(s, E)$.  It reads
\begin{align*}
\sum_{\rho} \phi \left( \gamma \frac{\log{X}}{2 \pi} \right) 
& =  \widehat{\phi}(0) \frac{\log{N}}{\log{X}} 
+ \frac{2}{2 \pi } \intR \frac{\gamma'}{\gamma}\left(\onehalf + it, E \right) \phi\left(t \frac{\log{X}}{2 \pi}\right) dt 
\\
& \qquad - \frac{2}{\log{X}} \sum_{n} \frac{\Lambda_E(n)}{n^{1/2}} \widehat{\phi} \left(\frac{\log{n}}{\log{X}} \right),
\end{align*}
where $\gamma(s, E) = (2 \pi)^{-s} \Gamma(s + \onehalf)$ and
\begin{equation*}
 \frac{L'}{L}(s, E) = -\sum_{n} \Lambda_E(n) n^{-s}.
\end{equation*}
The contribution from the integral of $\gamma'/\gamma$ is
\begin{equation}
\frac{2}{\log{X}} \int_{- \infty}^{\infty} \phi(x) \left\{ -\log{2\pi} + \frac{\Gamma'}{\Gamma} \left(1 + \frac{2 \pi i x}{\log{X}} \right) \right\}dx.
\end{equation}
For $p | N$ the logarithmic derivative of the local factor of $L$ is
\begin{equation*}
\frac{L_p'}{L_p}(s, E) = - \sum_{\nu = 1}^{\infty} \frac{\lambda^{\nu}(p)}{p^{\nu s}} \log{p} 
\end{equation*}
so the contribution from such a $p$ is
\begin{eqnarray*}
 -2 \frac{\log{p}}{\log{X}} \sum_{\nu = 1}^{\infty} \lambda^{\nu}(p) p^{-\nu/2} \; \widehat{\phi}\left( \nu \frac{\log{p}}{\log{X}}\right).
\end{eqnarray*}
The contribution from the primes not dividing $N$ can be expressed in various ways (e.g. $\lambda(p^\nu)$, $\lambda^\nu(p)$, or the parameters $\alpha(p)$, $\overline{\alpha(p)}$ such that $\alpha(p) + \overline{\alpha(p)} =  \lambda(p)$).  We choose to make the derivation with respect to $\lambda(p^\nu)$.  Since for $p \notdivtext N$ 
\begin{equation*}
L_p(s, E) = \left(1 - \frac{\lambda(p)}{p^{s}} + \frac{1}{p^{2s}} \right)^{-1} = \sum_{\nu = 0}^{\infty} \frac{\lambda(p^\nu)}{p^{\nu s}},
\end{equation*}
we therefore have
\begin{equation*}
\frac{L_p'}{L_p}(s, E) = - \log{p} \left(\frac{\lambda(p)}{p^s} - \frac{2}{p^{2s}}\right) \sum_{\nu = 0}^{\infty} \frac{\lambda(p^\nu)}{p^{\nu s}}.
\end{equation*}
Using $\lambda(p) \lambda(p^\nu) = \lambda(p^{\nu + 1}) + \lambda(p^{\nu - 1})$ (with the convention $\lambda(p^{-1}) = 0$) gives
\begin{align*}
\frac{L_p'}{L_p}(s, E) & = - \log{p} \left(\sum_{\nu = 0}^{\infty} \frac{\lambda(p^{\nu + 1})}{p^{(\nu + 1)s}} + \sum_{\nu = 0}^{\infty} \frac{\lambda(p^{\nu - 1})}{p^{(\nu + 1)s}} - 2 \sum_{\nu = 0}^{\infty} \frac{\lambda(p^{\nu})}{p^{(\nu + 2)s}} \right) \\
& = -\log{p} \left(\sum_{\nu = 1}^{\infty} \frac{\lambda(p^{\nu })}{p^{\nu s}} - \sum_{\nu = 0}^{\infty} \frac{\lambda(p^{\nu})}{p^{(\nu + 2)s}} \right)
\end{align*}
\begin{equation*}
= \frac{\log{p}}{p^{2s}} - \log{p} \sum_{\nu = 1}^{\infty} \lambda(p^\nu) \left(\frac{1}{p^{\nu s}} - \frac{1}{p^{(\nu + 2)s}} \right).
\end{equation*}
Therefore the local integral contribution from this $p$ is
\begin{equation*}
\frac{2\log{p}}{p \log{X}} \widehat{\phi}\left(  \frac{2 \log{p}}{\log{X}} \right) - 2 \frac{  \log{p}}{\log{X}}  \sum_{\nu = 1}^{\infty} \frac{\lambda(p^\nu)}{p^{\nu/2}} \left( \widehat{\phi}\left(  \frac{\nu \log{p}}{\log{X}} \right) - p^{-1} \widehat{\phi}\left(  \frac{(\nu + 2) \log{p}}{\log{X}} \right) \right).
\end{equation*}
Now we have proved
\begin{myprop}[Explicit formula for an elliptic curve]
\label{prop:explicitformula}
\begin{align*}
\sum_{\gamma} \phi\left(\gamma \frac{\log X}{2\pi} \right) 
& = \widehat{\phi}(0) \frac{\log{N}}{\log{X}} 
+ \frac{2}{\log{X}} \int_{- \infty}^{\infty} \phi(x) \left\{ -\log{2\pi} 
+ \frac{\Gamma'}{\Gamma} \left(1 + \frac{2 \pi i x}{\log{X}} \right) \right\}dx \\
& \qquad -  2 \sum_{p | N} \frac{\log{p}}{\log{X}} \sum_{\nu = 1}^{\infty} \frac{\lambda^{\nu}(p)}{p^{\nu/2}} \; \widehat{\phi}\left(  \frac{\nu \log{p}}{\log{X}}\right) 
+ \sum_{p \notdiv N} \frac{2\log{p}}{p \log{X}} \widehat{\phi}\left(  \frac{2 \log{p}}{\log{X}} \right) \\
& \qquad - 2 \sum_{p \notdiv N}  \frac{  \log{p}}{\log{X}}  \sum_{\nu = 1}^{\infty} \frac{\lambda(p^\nu)}{p^{\nu/2}} \left( \widehat{\phi}\left(  \frac{\nu \log{p}}{\log{X}} \right) - p^{-1} \widehat{\phi}\left(  \frac{(\nu + 2) \log{p}}{\log{X}} \right) \right).
\end{align*}
\end{myprop}

\section{Computing the Conductor and Dirichlet Coefficients}
\label{section:invariants}
In this section we calculate the arithmetical quantities we will need in order to apply the explicit formula.  

We briefly review some basic facts of elliptic curves.  Suppose the curve $E/\mq$ is given by the general Weierstrass equation
\begin{equation*}
E: y^2 + a_1 xy + a_3y = x^3 + a_2 x^2 + a_4 x + a_6.
\end{equation*}
We have the canonical parameters
\begin{gather*}
b_2 = a_1^2 + 4 a_2, \qquad
b_4 = 2a_4 + a_1 a_3, \qquad
b_6  = a_3^2 + 4a_6, \\
b_8  = a_1^2 a_6 + 4 a_2 a_6  - a_1 a_3 a_4 + a_2 a_3^2 - a_4^2, \\
c_4 = b_2^2 - 24b_4, \qquad
c_6 = -b_2^3 + 36b_2 b_4 - 216 b_6, \\
\Delta = -b_2^2 b_8 - 8 b_4^3 - 27 b_6^2 + 9 b_2 b_4 b_6.
\end{gather*}
The conductor $N$ of $E$ is a positive integer dividing $\Delta$.  There is no simple formula for $N$, but it can be computed by following an algorithm due to Tate (see \cite{SilvermanAEC2}, pp. 361-368 for a description).  If the Weierstrass equation defining $E$ is minimal then every prime dividing $\Delta$ also divides $N$ (recall that a Weierstrass equation is minimal if $\Delta$ cannot be reduced by a change of variables).  All curves in our families $\mathcal{F}_1$ and $\mathcal{F}_2$ are minimal.  The following formula determines the power of $p$ dividing $N$
\begin{equation*}
\text{ord}_p(N) =
\begin{cases}
0 & \text{if } (p, \Delta) = 1 \\
1 & \text{if } p | \Delta \text{ and } (p, c_4) = 1 \\
2 + \delta_p(E) & \text{if } p | (\Delta, c_4),
\end{cases}
\end{equation*}
where $\delta_p(E)$ is a measure of {\it wild ramification} that may occur at the primes $p =2$ and $3$.  If $p > 3$ then $\delta_p(E) = 0$.  The three cases above correspond to $E$ having good reduction at $p$, multiplicative reduction at $p$, and additive reduction at $p$, respectively.

For convenience of notation we make the following
\begin{mydefi}
Let $r_1$ be an integer such that the congruence $a \equiv r_1 \pmod{2q}$ is equivalent to the congruences $a \equiv a_0 \pmod{q}$ and $a \equiv 1 \pmod{2}$.  Similarly define $t_1$ with respect to $b$ and $b_0$.   Analogously, let $r_2$ and $t_2$ be integers such that the congruences $a \equiv a_0 \pmod{q}$ and $a \equiv 1 \pmod{4}$ are equivalent to $a \equiv r_2 \pmod{4q}$, and similarly for $b$.  
\end{mydefi}
Note that $r_1$ and $t_1$ codify the congruence conditions of the family $\mathcal{F}_1$ while $r_2$ and $t_2$ perform the same function for $\mathcal{F}_2$.  We may unify the conditions for the two families by noting that the conditions are special cases of the congruences $a \equiv r_i \pmod{2^i q}$, $b \equiv t_i \pmod{2^i q}$.

\subsection{The family $y^2 = x(x-a)(x+b)$}
Consider the elliptic curve $E$ given by the Weierstrass equation $y^2 = x(x-a)(x + b)$.  The curve $E$ has discriminant $\Delta = 16a^2 b^2(a + b)^2$  and parameter $c_4 = 16(a^2 + ab + b^2) = 16((a + b)^2 - ab)$.  If $a$ and $b$ are relatively prime then $(c_4, \Delta)$ is a power of $2$, and hence $E$ does not have additive reduction at any primes $p \neq 2$.
\begin{mylemma} 
\label{lem:conductor}
Suppose $a$ and $b$ are odd and coprime.  For an integer $n$ let $n^*$ be the product of primes dividing $n$. Then the curve $E_{a, b}: y^2 = x(x-a)(x + 2b)$ has conductor
\begin{equation*}
N = 2^5 a^* b^* (a + 2b)^*.
\end{equation*}
\end{mylemma}
Proof.  By the discussion before the lemma, the odd part of $N$ is simply the product of odd primes dividing $ab(a + 2b)$, which is $a^* b^* (a+2b)^*$ because $a$ and $b$ are odd and coprime. It remains to compute ord$_2(N)$.  We follow Tate's algorithm \cite{SilvermanAEC2}.  The point $(0, 0)$ is already singular $\pmod{2}$.  The algorithm terminates at step 4, since $2^3 \negmedspace \notdivtext b_8$ (note $b_8 = 4a^2 b^2$ and $a$ and $b$ are odd).  Tate's algorithm says ord$_2(N) = \text{ord}_2(\Delta) - 1$, which is $5$ since $2^{-6} \Delta = a^2 b^2 (a + 2b)^2$ is odd. \eop

The following lemma provides a useful formula for the Dirichlet coefficients of $L(s, E_{a, b})$.
\begin{mylemma}
\label{lem:lambdaformula}
Let $a, b, E_{a, b}$ be as in Lemma \ref{lem:conductor}.  Let $\lambda_{a, b}(p)$ be the coefficient of $p^{-s}$ in the Dirichlet series expansion of $L(s, E_{a, b})$.  If $p \neq 2$ then
\begin{equation}
\label{eq:lambdadef}
\lambda_{a, b}(p) =  - p^{-1/2} \sum_{x \shortmod{p}} \left(\frac{x(x-a)(x+2b)}{p} \right).
\end{equation}
If $p | N$ then 
\begin{equation}
\label{eq:lambdaformular}
p^{1/2} \lambda_{a, b}(p) = \begin{cases}
\left(\frac{2b}{p} \right) & \text{if } p | a \\
\left(\frac{-a}{p} \right) & \text{if } p | b \\
\left(\frac{2b}{p} \right)  = \left(\frac{-a}{p} \right) & \text{if } p | (a +2b) \\
0 & \text{if } p = 2.
\end{cases}
\end{equation}
\end{mylemma}
Proof.  The conditions $(a, b) = 1$, $a$ and $b$ odd imply that the Weierstrass equation defining $E_{a, b}$ is minimal, and hence that (\ref{eq:lambdadef}) holds.  We easily evaluate the sum when $p | N$.  \eop

\subsection{The family $y^2 = x(x^2 +ax - b)$}
Consider the elliptic curve $E$ given by the Weierstrass equation $y^2 = x(x^2 + ax - b)$.  $E$ has discriminant $\Delta = 16b^2(a^2 + 4b)$ and parameter $c_4 = 16(a^3 + 3b)$.  If $a$ and $b$ are relatively prime then $(c_4, \Delta)$ is a power of $2$.
\begin{mylemma}
\label{lem:conductor2}
Suppose $a$ and $b$ are coprime, and that $a \equiv b \equiv 1 \pmod{4}$.  Then the curve $E_{a, b}: y^2 = x(x^2 + 2ax - b)$ has conductor
\begin{equation*}
N = 2^{7} b^{*} \left(\frac{a^2 + b}{2}\right)^{*} = 2^6 b^* (a^2 + b)^*.
\end{equation*}
\end{mylemma}
Proof.  The odd part of $N$ is easily shown to be as claimed.  We need only compute ord$_2(N)$.  We follow Tate's algorithm.  First apply the change of variables $x \rightarrow x + 1$, giving the Weierstrass equation $E':y^2 = x^3 + x^2(2a + 3) + x(4a - b + 3) + (2a - b + 1)$.  The reduced curve $E' \pmod{2}$ then has singular point $(0, 0)$.  We have $a_6 = 2a - b + 1 \equiv  2 \pmod{4}$, so the algorithm terminates at Step 3.  Tate's algorithm reads that $\text{ord}_2(N) =  \text{ord}_2(\Delta)$. We have $\text{ord}_2(\Delta) = 7$ since $2^{-7} \Delta = b^2 ((a^2 + b)/2) \equiv 1 \pmod{2}$. \eop

\begin{mylemma}
\label{lem:lambdaformula2}
Let $a, b$, and $E_{a, b}$ be as in Lemma \ref{lem:conductor2}.  Let $\lambda_{a, b}(p)$ be the coefficient of $p^{-s}$ in the Dirichlet series expansion of $L(s, E_{a, b})$.  If $p \neq 2$ then
\begin{equation}
\lambda_{a, b}(p) =  - p^{-1/2} \sum_{x \shortmod{p}} \left(\frac{x(x^2 +2ax - b)}{p} \right).
\end{equation}
If $p | N$ then 
\begin{equation}
p^{1/2} \lambda_{a, b}(p) = \begin{cases}
 \left(\frac{2a}{p} \right) & \text{if } p | b \\
 \left(\frac{-a}{p} \right) & \text{if } p | (a^2 + b), p \neq 2 \\
0 & \text{if } p = 2 
\end{cases}
\end{equation}
\end{mylemma}
Proof.  Identical to that of Lemma \ref{lem:lambdaformula}. \eop

\section{Computing the Size of the Families}
\label{section:sizes}
In this section we compute the sizes of the families $\mathcal{F}_1$ and $\mathcal{F}_2$.  Essentially we need to count pairs of coprime integers $a$ and $b$ that lie in specified arithmetic progressions.
\subsection{A General Computation}
It is useful to make the following
\begin{mydefi}
For any integer $l \neq 0$ set
\begin{equation*} 
\gamma(l) = \frac{1}{l} \prod_{p | l} \frac{1}{1 - p^{-2}},
\end{equation*}
where the product is over the prime divisors of $l$.
\end{mydefi}
Note that $\gamma(l)$ is multiplicative.

\begin{mylemma} 
\label{lem:sizelemma}
Let $A$ and $B$ be larger than $2$. Let $r$, $t$, $l_1$, and $l_2$ be positive integers such that  $(r, l_1) = (t, l_2) = 1$.  Let $w \in C_0^{\infty}(\mr^2)$, $\widehat{w}(0,0) = 1$.  Then
\begin{equation*}
\mathop{\mathop{\mathop{\sum \sum}_{a \equiv r \shortmod{l_1}} }_{b \equiv t \shortmod{l_2}} }_{(a, b) = 1} w \left( \frac{a}{A}, \frac{b}{B} \right) = AB  \gamma(l_1 l_2) \zeta^{-1}(2) + O((A + B)^{1 + \varepsilon}).
\end{equation*}
Further, suppose 
$p$ is prime such that $(p, l_1 l_2) = 1$.  Then
\begin{equation*}
\mathop{\mathop{\mathop{ \mathop{\sum \sum}_{a \equiv 0 \shortmod{p}}}_{a \equiv r \shortmod{l_1}} }_{b \equiv t \shortmod{l_2}} }_{(a, b) = 1} w \left( \frac{a}{A}, \frac{b}{B} \right) = AB  \gamma(l_1 l_2) \zeta^{-1}(2) \frac{1}{p + 1} + O((A + B)^{1 + \varepsilon}).
\end{equation*}
The implied constants depend only on $w$ and $\varepsilon$.
\end{mylemma}
Note.  The second statement will not be applied until Section \ref{section:primesdividing} but it is natural to prove it here.

Proof.  We prove both statements simultaneously.  Let $q$ be either $1$ or $p$ and let $S$ be the sum to be calculated.
We have
\begin{equation*}
S = \mathop{\mathop{\mathop{ \mathop{\sum \sum}_{a \equiv 0 \shortmod{q}}}_{a \equiv r \shortmod{l_1}} }_{b \equiv t \shortmod{l_2}} }_{(a, b) = 1} w \left( \frac{a}{A}, \frac{b}{B} \right)
=
\mathop{\mathop{\mathop{\sum \sum}_{a \equiv \bar{q} r \shortmod{l_1}} }_{b \equiv t \shortmod{l_2}} }_{(aq, b) = 1} w \left( \frac{aq}{A}, \frac{b}{B} \right).
\end{equation*}
By M\"{o}bius inversion, we obtain
\begin{align*}
 S & = \mathop{\sum_{(d, l_1 l_2 q) = 1}}_{d \ll \text{min}(A, B)} \mu(d) \mathop{ \mathop{\mathop{\sum \sum}_{a \equiv \overline{dq} r \shortmod{l_1}} }_{b \equiv \overline{d} t \shortmod{l_2}} }_{(b, q) = 1} w \left( \frac{adq}{A}, \frac{bd}{B} \right) \\
& = \mathop{\sum_{(d, l_1 l_2 q) = 1}}_{d \ll \text{min}(A, B)} \mu(d) \sum_{e | q} \mu(e) \mathop{\mathop{\sum \sum}_{a \equiv \overline{dq} r \shortmod{l_1}} }_{b \equiv \overline{de} t \shortmod{l_2}} w \left( \frac{adq}{A}, \frac{bde}{B} \right).
\end{align*}
Completing the sum in $a \pmod{l_1}$ and $b \pmod{l_2}$ gives
\begin{equation*}
S = \frac{AB}{l_1 l_2 q} \mathop{\sum_{d \ll \text{min}(A, B)}}_{(d, l_1 l_2 q) = 1} \frac{\mu(d)}{d^2} \sum_{e | q} \frac{\mu(e)}{e} \sum_h \sum_k  e\left( \frac{\overline{dq} r h}{l_1} \right) e\left( \frac{ \overline{de}t k}{l_2} \right) \widehat{w}\left(\frac{hA}{d ql_1}, \frac{kB}{del_2} \right).
\end{equation*}
Extracting the main term from the zero frequencies gives
\begin{align*}
S & = \frac{AB \varphi(q)}{l_1 l_2 q^2} \mathop{\sum_{d \ll \text{min}(A, B)}}_{(d, l_1 l_2 q) = 1} \frac{\mu(d)}{d^2} \widehat{w}\left(0, 0 \right) + O((A + B)^{1 + \varepsilon}) \\
& =  AB  \gamma(l_1 l_2) \gamma(q^2) \varphi(q) \zeta^{-1}(2) + O((A + B)^{1 + \varepsilon}), 
\end{align*}
since
\begin{equation*}
\mathop{\sum_{d=1}^{\infty}}_{(d, p_1 \ldots p_k)=1} \frac{\mu(d)}{d^s} = \frac{1}{\zeta(s)} \prod_{i = 1}^{k} \left(\frac{1}{1 - p_i^{-s}} \right)
\end{equation*}
if $p_1, \ldots, p_k$ are distinct primes.  We easily compute $\varphi(q) \gamma(q^2) = (p + 1)^{-1}$ in case $q = p$, and $\varphi(q) \gamma(q^2) = 1$ in case $q = 1$. \eop

\subsection{The Family $\mathcal{F}_1$}
\begin{mylemma}
\label{lem:mass}
The size of the family $\mathcal{F}_1$ is given by
\begin{equation*}
\sum_{E_{a, b} \in \mathcal{F}_1} w_X(E_{a, b}) = \frac{AB  \gamma(q^2)}{3 \zeta(2)} + O\left(X^{1/3 + \varepsilon} \right),
\end{equation*}
the implied constant depending only on $w$ and $\varepsilon$.
\end{mylemma}
The proof is a simple application of Lemma \ref{lem:sizelemma}.  By definition,
\begin{equation*}
\sum_{E_{a, b} \in \mathcal{F}_1} w_X(E_{a, b})  = \mathop{ \mathop{\mathop{\sum \sum}_{a \equiv r_1 \shortmod{2q}} }_{b \equiv t_1 \shortmod{2q}} }_{(a, b) = 1} w \left( \frac{a}{A}, \frac{b}{B} \right).
\end{equation*}
Simply apply Lemma \ref{lem:sizelemma} with $l_1 = l_2 = 2q$.  We easily compute $\gamma(4q^2) = \gamma(q^2)/3$. \eop

\subsection{The Family $\mathcal{F}_2$}
\begin{mylemma}
\label{lem:mass2}
The size of the family $\mathcal{F}_2$ is given by
\begin{equation*}
\sum_{E_{a, b} \in \mathcal{F}_2} w_X(E_{a, b}) = \frac{AB  \gamma(q^2)}{12 \zeta(2)}+ O\left(X^{1/2 + \varepsilon} \right),
\end{equation*}
the implied constant depending only on $w$ and $\varepsilon$.
\end{mylemma}
Proof.  By definition,
\begin{equation*}
\sum_{E_{a, b} \in \mathcal{F}_2}  w_X(E_{a, b}) = \mathop{ \mathop{\mathop{\sum \sum}_{a \equiv r_2 \shortmod{4q}} }_{b \equiv t_2 \shortmod{4q}} }_{(a, b) = 1} w \left( \frac{a}{A}, \frac{b}{B} \right).
\end{equation*}
We compute $\gamma(16q^2) = \gamma(q^2)/12$. \eop 

We now make the following
\begin{mydefi}  Set
\begin{eqnarray*}
M(\mathcal{F}_1) = \frac{AB  \gamma(q^2)}{3 \zeta(2)} = \frac{AB  \gamma(q)}{3\zeta(2) q}, \\
M(\mathcal{F}_2) = \frac{AB  \gamma(q^2)}{12 \zeta(2)} = \frac{AB  \gamma(q)}{12 \zeta(2)q}.
\end{eqnarray*}
\end{mydefi}
Of course the prevous two lemmas showed $W_X(\mathcal{F}_i) \sim M(\mathcal{F}_i)$.

\section{The Conductor on Average}
In this section we average the logarithm of the conductors of the families $\mathcal{F}_1$ and $\mathcal{F}_2$.  This is the technically most difficult aspect of this paper.  We are able to succeed with these two families essentially because the conductor factors into linear polynomials (roughly speaking, $N \approx ab(a + 2b)$).  The analogous calculation for the family of all elliptic curves would involve summing the logarithm of the product of primes dividing $4a^3 + 27 b^2$ (the summation being over $a$ and $b$ in appropriate ranges).  Getting a precise asymptotic for this sum is an open problem.

\label{section:conductor}
\subsection{A Sum with Squarefree Integers}
In this section we evaluate a sum involving the squarefree part of integers.  The calculation will be applied many times with various test functions and congruence conditions.  We have
\begin{mylemma}
\label{lem:squarefreecalculation}
Let $A$ and $B$ be larger than 2.  Let $r$ and $t$ be integers, and suppose $l$ is a positive integer such that  $(r, l) = (t, l) = 1$.  Let $w \in C_0^{\infty}(\mr^+ \times \mr^+)$, $\widehat{w}(0, 0) =1$.  Let $n^*$ be the product of primes dividing $n$.  Then
\begin{align*}
 \mathop{ \mathop{\mathop{\sum \sum }_{a \equiv r \shortmod{l}} }_{b \equiv t \shortmod{l}} }_{(a, b) = 1} \log{a^*} \, w \left( \frac{a}{A}, \frac{b}{B} \right) 
=  \frac{AB \gamma(l)}{\zeta(2) l} \left( \log{A} + 
\int \int \log{x} \, w(x, y) dx dy - \sum_{p \notdiv l} \frac{\log{p}}{p^2 - 1} \right) \\
 + \; O((A + A^{1/2}B)A^{\varepsilon}).
\end{align*}
The implied constant depends on $w$ and $\varepsilon$ only.
\end{mylemma}
By Lemma \ref{lem:sizelemma} notice that $AB \gamma(l) l^{-1} \zeta^{-1}(2)$ is the asymptotic if $\log{a^*}$ is not present.  The expected extra factor of $\log{A}$ to compensate for $\log{a^*}$ is correct but there is a rather large secondary term present.  

The outline of the proof is to execute the summation over $b$ and reduce to a one-variable sum, to which we will apply Lemma \ref{lem:onevar}.  

Proof.  Let $T$ be the sum to be computed.  By M\"{o}bius inversion we have
\begin{equation*}
T = \mathop{ \mathop{\mathop{\sum \sum }_{a \equiv r \shortmod{l}} }_{b \equiv t \shortmod{l}} }_{(a, b) = 1} \log{a^*} \, w \left( \frac{a}{A}, \frac{b}{B} \right) 
= \sum_{a \equiv r \shortmod{l}} \log{a^*} \sum_{d | a}  \mu(d) \sum_{b \equiv \overline{d} t \shortmod{l}} w \left( \frac{a}{A}, \frac{bd}{B} \right).
\end{equation*}
Applying Poisson summation in $b \pmod{l}$ gives
\begin{equation*}
T = \mathop{\sum_{a \equiv r \shortmod{l}}}_{a \ll A} \log{a^*} \sum_{d | a} \mu(d) \left( \frac{B}{dl} W \left( \frac{a}{A} \right)  + O(1) \right),
\end{equation*}
where
\begin{equation*}
W(x) = \int_{-\infty}^{\infty} w(x, y) dy.
\end{equation*}
Evaluating the summation over $d$ gives
\begin{equation*}
T = \frac{B}{l} \sum_{a \equiv r \shortmod{l}} \frac{\varphi(a)}{a} \log{a^*} \, W \left( \frac{a}{A} \right) + O(A^{1 + \varepsilon}).
\end{equation*}
The following lemma allows us to execute the summation over $a$ and will complete the proof.  \eop
\begin{mylemma}
\label{lem:onevar}
Let $A$ be larger than 2, let $r$ be an integer, and suppose $l$ is a positive integer such that  $(r, l) = 1$.  Let $W \in C_0^{\infty}(\mr^+)$, $\widehat{W}(0) = 1$.  Set
\begin{equation*}
S = \sum_{a \equiv r \shortmod{l}} \frac{\varphi(a)}{a} \log{a^*} \, W\left( \frac{a}{A}\right). 
\end{equation*}
Then
\begin{equation*}
S = \frac{A \gamma(l)}{\zeta(2)} \left( \log{A}+  \int_{0}^{\infty} \log{x} \, W(x) dx - \sum_{p \notdiv l} \frac{\log{p}}{p^2 - 1}  \right) + O(A^{1/2 + \varepsilon}),
\end{equation*}
the implied constant depending on $w$ and $\varepsilon$.
\end{mylemma}

Proof.  The overall idea is to use zeta function theory to evaluate $S$ (that is, we use the Mellin transform of $W$ to relate $S$ to the values of a certain zeta function $Z(s)$ and its logarithmic derivative at the point $s=1$ (see \eqref{eq:zetaexpression}).  

Write $a = a_1 a_2$ where $a_1$ is squarefree and divisible by $a_2^*$ (such a representation is unique since the conditions imply $a_1 = a^*$).  We obtain
\begin{equation*}
S =  \sum_{(a_2, l)=1} \mathop{\sum_{a_1 \equiv r \overline{a_2} \shortmod{l}}}_{a_1 \equiv 0 \shortmod{a_2^*}} \mu^2 (a_1) \frac{\varphi(a_1)}{a_1} \log{a_1}  \; W\left(\frac{a_1 a_2}{A} \right).
\end{equation*}
Now write $a_1 = a_2^* a_3$ (hence $a_2^*$ and $a_3$ are coprime and both squarefree) and set $\mathop{u \equiv r \overline{a_2 a_2^*} \pmod{l}}$ to obtain
\begin{equation}
\label{eq:Ssum}
S =  \sum_{(a_2, l) = 1}  \frac{\varphi(a_2)}{a_2} \mathop{\sum_{a_3 \equiv u \shortmod{l}}}_{(a_3, a_2^*) = 1} \mu^2 (a_3) \frac{\varphi(a_3)}{a_3} \log{(a_2^* a_3)}  W\left(\frac{a_2 a_2^* a_3}{A} \right).
\end{equation}
Let $S_{1}$ be the inner sum over $a_3$ in (\ref{eq:Ssum}) and set
\begin{equation*}
f(x) = \log (a_2^* x) W\left( \frac{a_2 a_2^* x}{A} \right).
\end{equation*}
By Mellin inversion,
\begin{equation*}
f(x) = \frac{1}{2 \pi i} \int_{(\sigma)} \left(\frac{A}{a_2 a_2^*}\right)^s F(s) x^{-s} ds,
\end{equation*}
where
\begin{equation*}
F(s) = \int_{0}^{\infty} \log\left(\frac{xA}{a_2}\right) W(x) x^{s-1} dx.
\end{equation*}
Then
\begin{align*}
S_{1} & = \mathop{\sum_{a_3 \equiv u \shortmod{l}}}_{(a_3, a_2^*) = 1} \mu^2 (a_3)  \frac{\varphi(a_3)}{a_3} f(a_3) \\
& = \frac{1}{2\pi i} \int_{(\sigma)}  \left(\frac{A}{a_2 a_2^*}\right)^s F(s) Z_{1}(s) ds,
\end{align*}
where
\begin{equation*}
Z_{1} (s) = \mathop{\sum_{n \equiv u \shortmod{l}}}_{(n, a_2^*) = 1} \mu^2 (n) \frac{\varphi(n)}{n} n^{-s}.
\end{equation*}
Unfortunately, $Z_{1}$ does not have an Euler product because of the restriction $n \equiv u \pmod{l}$ so write
\begin{eqnarray*}
Z_{1}(s) & = & \frac{1}{\varphi(l)} \sum_{\chi \shortmod{l}} \overline{\chi}(u) \sum_{(n, a_2^*) = 1} \mu^2 (n) \chi(n) \frac{\varphi(n)}{n} n^{-s} \\
& = & \frac{1}{\varphi(l)} \sum_{\chi \shortmod{l}} \overline{\chi}(u) Z_{\chi}(s),
\end{eqnarray*}
say.  $Z_\chi(s)$ has the following Euler product expansion
\begin{eqnarray*}
Z_\chi(s) &=& \prod_{p \notdiv a_2} \left(1 + \chi(p) (1 - p^{-1}) p^{-s} \right) \\
&=& \prod_{p} \left(1 + \chi(p) (1 - p^{-1}) p^{-s} \right) \prod_{p | a_2} \left(1 + \chi(p) (1 - p^{-1}) p^{-s} \right)^{-1}.
\end{eqnarray*}
The product
\begin{equation*}
V_{\chi}(s) := \prod_{p} \left(1 + \chi(p) (1 - p^{-1}) p^{-s} \right) \left( 1 - \chi(p) p^{-s}\right) = \prod_p \left( 1 - \frac{\chi(p)}{p^{s+1}} - \frac{\chi(p^2)}{p^{2s}} + \frac{\chi(p^2)}{p^{2s + 1}} \right)
\end{equation*}
is absolutely convergent for Re $s > 1/2$ so that
\begin{equation*}
Z_{\chi}(s) = L(s, \chi) V_{\chi}(s) \prod_{p | a_2} \left(1 + \chi(p) (1 - p^{-1}) p^{-s} \right)^{-1}
\end{equation*}
is holomorphic in Re $s > 1/2$ except for a simple pole at $s=1$ that occurs if and only if $\chi$ is the principal character $\psi_l \pmod{l}$.  We evaluate $S_{1}$ by moving the contour of integration to the line Re $s = 1/2 + \varepsilon$, picking up a pole at $s=1$ for $Z_{\psi_l}(s)$.  On the line $\text{Re }s = 1/2 + \varepsilon$ we have the bounds
\begin{align*}
L(s, \chi) & \ll |s| |l|^{3/16 + \varepsilon}, \\
V_{\chi}(s) & \ll 1, \\
 \text{ and } & \\
\prod_{p | a_2} \left(1 + \chi(p) (1 - p^{-1}) p^{-s} \right)^{-1} & \ll |a_2|^{\varepsilon/4}
\end{align*}
uniformly in $l$ (using the Burdgess bound for $L(s, \chi)$), and of course
\begin{equation*}
F(s) \ll |s|^{-100} (1 + (a_2 A)^{\varepsilon/4}).
\end{equation*}
Thus we obtain
\begin{eqnarray*}
S_{1} & = & \frac{1}{\varphi(l)} \sum_{\chi \shortmod{l}} \overline{\chi}(u) \frac{1}{2 \pi i} \int_{(\sigma)} \left(\frac{A}{a_2 a_2^*}\right)^s F(s) Z_{\chi}(s) ds \\
& = & \frac{1}{\varphi(l)} \frac{A}{a_2 a_2^*} \psi_l(u) F(1) \, \text{Res}_{s = 1} Z_{\psi_l}(s) + O \left( \left(\frac{A}{a_2 a_2^*}\right)^{1/2 + \varepsilon} (a_2 A)^{\varepsilon/2} \right),
\end{eqnarray*}
uniformly in $l$.
Inserting this expression for $S_1$ into (\ref{eq:Ssum}), we get
\begin{equation}
S = \frac{A}{ \varphi(l)} \sum_{(a_2, l) = 1} \frac{\varphi(a_2)}{a_2} \frac{1}{a_2 a_2^*} F(1) \, \text{Res}_{s=1} Z_{\psi_l}(s) + O(A^{1/2 + \varepsilon}).
\end{equation}
Setting $V_{l}(s) = V_{\psi_l}(s)$ we obtain
\begin{equation*}
\text{Res}_{s=1} Z_{\psi_l}(s) = \frac{\varphi(l)}{l} V_l(1) \prod_{p | a_2} \left(1 + (1 - p^{-1}) p^{-1} \right)^{-1},
\end{equation*}
since $\text{Res}_{s=1} L(s, \psi_l) = l^{-1} \varphi(l)$ and $\psi_l(p) = 1$ for all $p |  a_2$ (because $(a_2, l) = 1$).
Also,
\begin{equation*}
F(1) = \int_{0}^{\infty} \log\left( \frac{xA}{a_2}\right) W(x)dx. 
\end{equation*}
Hence
\begin{eqnarray*}
S = \frac{A}{l} V_l (1)  \sum_{(n, l) = 1}  \frac{\varphi(n)}{n} \frac{1}{n} \prod_{p | n} p^{-1} \left(1 + (1 - p^{-1}) p^{-1} \right)^{-1}  \int_0^{\infty}  \log\left( \frac{xA}{n}\right) W(x)dx \\
+ O(A^{1/2 + \varepsilon}).
\end{eqnarray*}
It is natural to define
\begin{equation*}
Z(s) = \sum_{(n, l)= 1}  \frac{\varphi(n)}{n} \left( \prod_{p | n} p^{-1} \left(1 + (1 - p^{-1}) p^{-1} \right)^{-1} \right) n^{-s},
\end{equation*}
for then
\begin{equation}
\label{eq:zetaexpression}
S = \frac{A}{l} V_l(1) Z(1) \left( \int_0^\infty \log(xA) W(x) dx + \frac{Z'}{Z}(1) \int_0^{\infty} W(x) dx \right) + O(A^{1/2 + \varepsilon}).
\end{equation}
A computation shows that $Z(s)$ has the Euler product expansion
\begin{eqnarray*}
Z(s) & = & \prod_{p \notdiv l} \left(1 + \frac{1}{p} \; \frac{1 - \frac{1}{p}}{1 +  \frac{1}{p} - \frac{1}{p^2}} \;  \frac{p^{-s}}{1 - p^{-s}} \right) \\
& = &  \prod_{p \notdiv l} \left(1 + \frac{1}{p} - \frac{1}{p^2} \right)^{-1} \left(1 + \frac{1}{p} (1 - \frac{1}{p}) (1 - p^{-s})^{-1} \right).
\end{eqnarray*}
We easily compute
\begin{eqnarray*}
V_l(1) Z(1) & =& \prod_{p \notdiv l} \left( 1 + \left(1 - \frac{1}{p} \right) \frac{1}{p} \right) \left( 1 - \frac{1}{p} \right) \prod_{p \notdiv l} \left[ \left(1 + \frac{1}{p} - \frac{1}{p^2} \right)^{-1} \left(1 + \frac{1}{p} \right) \right] \\
& = & \prod_{p \notdiv l} \left(1 - \frac{1}{p^2} \right) = \frac{l \gamma(l)}{\zeta(2)}.
\end{eqnarray*}
We compute
\begin{equation*}
\frac{Z'}{Z}(s) = - \sum_{p \notdiv l} \frac{\frac{1}{p}(1 - \frac{1}{p}) \log{p} }{1 + \frac{1}{p}(1 - \frac{1}{p}) (1 - p^{-s})^{-1}} (1 - p^{-s})^{-2} p^{-s},
\end{equation*}
so
\begin{equation*}
\frac{Z'}{Z}(1) = - \sum_{p \notdiv l} \frac{\log{p}}{p^2 - 1}.
\end{equation*}
We now conclude
\begin{equation*}
S = \frac{A \gamma(l)}{\zeta(2)} \left( \log{A}+ \int_{0}^{\infty} \log{x} W(x) dx - \sum_{p \notdiv l} \frac{\log{p}}{p^2 - 1}  \right) + O(A^{1/2 + \varepsilon}),
\end{equation*}
which completes the proof. \eop

\subsection{The Average Conductor of the Family $\mathcal{F}_1$}
In this section we evaluate the conductor on average.  We have
\begin{mylemma} 
\label{lem:conductorestimate}
Let $N_{a, b}$ be the conductor of the curve $y^2 = x(x-a)(x + 2b)$.  Then for $\varepsilon > 0$ we have
\begin{equation}
 \sum_{E_{a, b} \in \mathcal{F}_1} \frac{\log{N_{a, b}}}{\log{X}} w_X(E_{a, b}) =  \left\{1  + \frac{d_{1, 1}(q)}{\log{X}} + \frac{c_{1, 1}}{\log{X}}\right\} M(\mathcal{F}_1) + O(X^{1/2 + \varepsilon}),
\end{equation}
where
\begin{equation}
\label{eq:d1}
d_{1, 1}(q) = 3 \sum_{p | q} \frac{\log{p}}{p^2 - 1}
\end{equation}
and
\begin{equation}
\label{eq:c1}
c_{1, 1} =  \int_0^\infty \int_0^\infty \log(2^5 xy(x+2y)) w(x, y) dx dy - 3 \sum_{p \neq 2} \frac{\log{p}}{p^2 - 1}.
\end{equation}
The implied constant depends only on $w$ and $\varepsilon$.
\end{mylemma}
Note $d_{1, 1}(1) = 0$.

Proof.  
Using Lemma \ref{lem:conductor} and the additivity of the logarithm we write 
\begin{equation*}
\mathop{ \mathop{ \mathop{\sum \sum}_{a \equiv r_1 \shortmod{2q}} }_{b \equiv t_1 \shortmod{2q}} }_{(a, b) = 1} \log{(2^5 a^* b^* (a + 2b)^*)} w\left( \frac{a}{A}, \frac{b}{B} \right) =  S_{0} + S_{1} + S_{2} + S_{3},
\end{equation*}
where $S_0, S_1, S_2$, and $S_3$ correspond to the terms with the logarithm factor $\log{2^5}$, $\log{a^*}$, $\log{b^*}$, and $\log{(a + 2b)^*}$, respectively.
By Lemma \ref{lem:mass} we have
\begin{equation*}
S_0 =   \log{(2^5)} M(\mathcal{F}_1) + O(X^{1/3 + \varepsilon}).
\end{equation*}
An application of Lemma \ref{lem:squarefreecalculation} shows
\begin{equation*}
S_1 = M(\mathcal{F}_1) \left( \log{A} + \int_0^\infty \int_0^\infty \log{x} \; w(x, y) dx dy - \sum_{p \notdiv 2q} \frac{\log{p}}{p^2 - 1} \right)
+O(X^{1/2 + \varepsilon}).
\end{equation*}
Similarly, 
\begin{equation*}
S_2  =  M(\mathcal{F}_1) \left( \log{B} + \int_0^\infty \int_0^\infty \log{y} \; w(x, y) dx dy - \sum_{p \notdiv 2q} \frac{\log{p}}{p^2 - 1} \right) 
+ O(X^{1/2 + \varepsilon}).
\end{equation*}
The estimation for $S_3$ is only slightly different.  We first change variables via $a \rightarrow a - 2b$; this leads to (recall $A=B$)
\begin{equation*}
S_3 = \mathop{ \mathop{ \mathop{\sum \sum}_{a \equiv r_1 +2 t_1 \shortmod{2q}} }_{b \equiv t_1 \shortmod{2q}} }_{(a, b) = 1}  \log{a^*} \, w^*\left( \frac{a}{A}, \frac{b}{B} \right),
\end{equation*}
where $w^*(x, y) = w(x-2y, y)$.  We easily see $\widehat{w^*}(0, 0) = \widehat{w}(0, 0)$ and
\begin{equation*}
\int_0^{\infty} \log{x} \int_{0}^{\infty} w^*(x, y) dy dx = \int_0^\infty  \int_0^\infty \log(x + 2y) w(x, y) dx dy.
\end{equation*}
Now we've proved
\begin{equation*}
S_0 + S_1 + S_2 + S_3 = M(\mathcal{F}_1) \left( \log{X} +  d_{1, 1}(q)  +  c_{1, 1}\right) 
+ O(X^{1/2 +  \varepsilon})
\end{equation*}
where
\begin{equation*}
c_{1, 1} =  \int_0^\infty \int_0^\infty \log(2^5 xy(x+2y)) w(x, y) dx dy - 3 \sum_{p \neq 2} \frac{\log{p}}{p^2 - 1} 
\end{equation*}
is a constant depending on $w$ only.  Now the proof of Lemma \ref{lem:conductorestimate} is complete. \eop

\subsection{The Average Conductor of the Family $\mathcal{F}_2$}

We have
\begin{mylemma} 
\label{lem:conductorestimate2}
Let $N_{a, b}$ be the conductor of the curve $y^2 = x(x^2 + 2a x - b)$.  Then for $\varepsilon > 0$ we have
\begin{equation}
 \sum_{E_{a, b} \in \mathcal{F}_2} \frac{\log{N_{a, b}}}{\log{X}} w_X(E_{a, b}) =  \left\{1  + \frac{d_{1, 2} (q)}{\log{X}} + \frac{c_{1, 2}}{\log{X}} \right\} M(\mathcal{F}_2) +  O(X^{1/3 + \varepsilon}),
\end{equation}
where
\begin{equation}
\label{eq:d1'}
d_{1, 2}(q) = 2 \sum_{p | q} \frac{\log{p}}{p^2 - 1}
\end{equation}
and
\begin{equation}
\label{eq:c1'}
c_{1, 2} =  \int_0^\infty \int_0^\infty \log(2^6 y(x^2+y)) w(x, y) dx dy -  2 \sum_{p \neq 2} \frac{\log{p}}{p^2 - 1}.
\end{equation}
The implied constant depends only on $w$ and $\varepsilon$.
\end{mylemma}
Proof.  As in the proof of Lemma \ref{lem:conductorestimate}, break the sum up into $S_0 + S_1 + S_2$, where $S_0$ corresponds to $\log{2^6}$, $S_1$ corresponds to $\log{b^*}$, and $S_2$ corresponds to $\log{(a^2 + b)^*}$.  Obviously we have
\begin{equation*}
S_0 = \log{(2^6)} \, M(\mathcal{F}_2) + O(X^{1/2 + \varepsilon}).
\end{equation*}
By Lemma \ref{lem:squarefreecalculation} we easily have
\begin{equation*}
S_1 = M(\mathcal{F}_2) \left( \log{B} + \int_0^\infty \int_0^\infty \log{y} \; w(x, y) dx dy - \sum_{p \notdiv 2q} \frac{\log{p}}{p^2 - 1} \right)
+O(X^{1/2 + \varepsilon}).
\end{equation*}
To handle $S_2$ we apply the change of variables $b \rightarrow b - a^2$ and obtain
\begin{equation*}
S_2 = \mathop{ \mathop{ \mathop{ \sum \sum}_{a \equiv r_2 \shortmod{4q}} }_{b \equiv t_2 + r_2^2 \shortmod{4q}} }_{(a, b) = 1} \log{b^*} w^*\left( \frac{a}{A}, \frac{b}{B} \right),
\end{equation*}
where $w^*(x, y) = w(x, y -x^2)$.  An application of Lemma \ref{lem:squarefreecalculation} gives
\begin{equation*}
S_2 = M(\mathcal{F}_2) \left( \log{B} + \int_0^\infty \int_0^\infty \log{(x^2 + y)} \; w(x, y) dx dy - \sum_{p \notdiv 2q} \frac{\log{p}}{p^2 - 1} \right)
+O(X^{1/2 + \varepsilon}).
\end{equation*}
Adding $S_0$, $S_1$, and $S_2$ completes the proof. \eop

\section{The Gamma Factor}
The gamma factor in the functional equation of $L(s, E)$ is the same for every elliptic curve $E$ so there will be no variation in our families.
\begin{mylemma} For any integer $M \geq 1$ we have
\begin{equation*}
\frac{2}{\log{X}} \int_{-\infty}^\infty \phi(x) \left\{-\log{2\pi} + \frac{\Gamma'}{\Gamma} \left(1 + \frac{2\pi i x}{\log{X}}\right) \right\} dx  
\end{equation*}
\begin{equation*}
=  \frac{\widehat{\phi}(0)}{\log{X}}c_2 
- 2 \sum_{l = 1}^{M}  \zeta(1 + 2l)  \widehat{\phi}^{(2l)}(0) \left(\log{X} \right)^{-2l - 1} 
+ O_M\left((\log{X})^{-2M - 2}\right),
\end{equation*}
where
\begin{equation}
\label{eq:c2}
c_{2} = -2 \log{2\pi} - 2 \gamma := c_{2, 1} := c_{2, 2}.
\end{equation}
\end{mylemma}

Proof.  Set 
\begin{equation*}
I(X) = \int_{-\infty}^\infty \phi(x) \frac{\Gamma'}{\Gamma} \left(1 + \frac{2\pi i x}{\log{X}}\right) dx.
\end{equation*}
The notation $\psi = \frac{\Gamma'}{\Gamma}$ is standard.  We use the representations
\begin{equation*}
\psi(t) = \log{t} - \sum_{k = 0}^{\infty} \left[\frac{1}{t+k} - \log{\left(1 + \frac{1}{t+k}\right)} \right]
\end{equation*}
and
\begin{equation*}
\psi(1 + t) = -\gamma + \sum_{k = 2}^{\infty} (-1)^k \zeta(k) t^{k-1}, \; \; \; \; \; |t| < 1,
\end{equation*}
which are (8.362.2) and (8.363.1) in \cite{GR}, respectively.  The first clearly shows
\begin{equation*}
\psi \left(1 + \frac{2\pi i x}{\log{X}}\right) \ll \log{\left(1 + \frac{|x|}{\log{X}} \right)}.
\end{equation*}
Then
\begin{align*}
I(X) & = \int_{|x| \leq \frac{\log{X}}{4\pi}} + \int_{|x| > \frac{\log{X}}{4\pi}} \\
& = \int_{|x| \leq \frac{\log{X}}{4\pi}} \phi(x) \left( - \gamma - \sum_{l=1}^{\infty}  \zeta(1 + 2l) \left(\frac{2\pi ix}{\log{X}} \right)^{2l} \right) + O_Y\left( (\log{X})^{-Y}\right),
\end{align*}
where the odd powers of $x$ do not appear because $\phi$ is an even function.  Truncating the series at $l = N$ introduces an error of order $ (\log{X})^{-2N - 2}$.  Extending the integration back to the entire real line does not introduce a new error term.  Thus we obtain
\begin{equation*}
I(X) = - \gamma \widehat{\phi}(0) - \sum_{l = 1}^{N}  \zeta(1 + 2l) \left(\log{X} \right)^{-2l} \int_{-\infty}^{\infty} \phi(x) (2\pi ix)^{2l} dx + O_N\left((\log{X})^{-2N - 2}\right)
\end{equation*}
\begin{equation*}
= - \gamma \widehat{\phi}(0) - \sum_{l = 1}^{N}  \zeta(1 + 2l)  \widehat{\phi}^{(2l)}(0) \left(\log{X} \right)^{-2l} + O_N\left((\log{X})^{-2N - 2}\right),
\end{equation*}
since
\begin{equation*}
\int_{-\infty}^{\infty} \phi(x) (-2\pi ix)^{m} dx = \widehat{\phi}^{(m)}(0).
\end{equation*}
Now the proof is complete. \eop

\section{The Contribution From The Primes Dividing The Conductor}
\label{section:primesdividing}
\subsection{An Application of P\'{o}lya-Vinogradov}
It is convenient to state here the following
\begin{mylemma}
\label{lem:polya}
Let $A$ and $B$ be larger than $2$. Let $r$, $t$, $l_1$, and $l_2$ be positive integers such that  $(r, l_1) = (t, l_2) = 1$.  Let $w \in C_0^{\infty}(\mr^+ \times \mr^+)$, $\widehat{w}(0, 0) = 1$.  Let $\psi$ be a non-principal Dirichlet character of modulus $q$, where $q$ is coprime with $ l_2$.  Then
\begin{equation*}
\mathop{\mathop{\mathop{ \mathop{{\sum \sum}}_{a \equiv 0 \shortmod{q}}}_{a \equiv r \shortmod{l_1}} }_{b \equiv t \shortmod{l_2}} }_{(a, b) = 1} \psi(b) w \left( \frac{a}{A}, \frac{b}{B} \right) \ll l_2^{1/2} q^{-1/2} A \log{A},
\end{equation*}
the implied constant depending only on $w$.
\end{mylemma}
Proof.  Let $S$ be the sum to be estimated.  We may clearly assume $(q, l_1) =1$.  Then
\begin{align*}
S & = \mathop{\mathop{\mathop{ \sum \sum}_{a \equiv \overline{q} r \shortmod{l_1}} }_{b \equiv t \shortmod{l_2}} }_{(a, b) = 1} \psi(b) w \left( \frac{aq}{A}, \frac{b}{B} \right) \\
& = \mathop{\sum_{a \equiv \overline{q} r \shortmod{l_1}} }_{a \ll \frac{A}{q} } \sum_{d | a} \mu(d) \psi(d) \sum_{b \equiv \overline{d} t \shortmod{l_2}} \psi(b) w \left( \frac{aq}{A}, \frac{bd}{B} \right).
\end{align*}
By P\'{o}lya-Vinogradov, the summation over $b$ is $\ll (l_2 q)^{1/2}$.  Summing this bound over $a$ and $d$ completes the proof. \eop
\subsection{The Family $\mathcal{F}_1$}
\begin{mylemma}  
\label{lem:Tlemma}
Let $\lambda_{a, b}(p)$ be the $p$-th Dirichlet coefficient of $L(s, E_{a, b})$ for $E_{a, b} \in \mathcal{F}_1$.  Set
\begin{equation*}
T_1(a, b) = - \frac{2}{\log{X}} \sum_{p | N_{a, b}} \log{p} \sum_{\nu = 1}^{\infty} \lambda_{a, b}^{\nu}(p) p^{-\nu/2} \; \widehat{\phi}\left( \nu \frac{\log{p}}{\log{X}}\right) .
\end{equation*}
Then
\begin{equation}
\label{eq:Tsum}
\sum_{E_{a, b} \in \mathcal{F}_1} T_1(a, b) w_X(E_{a, b}) = M(\mathcal{F}_1) \frac{\widehat{\phi}(0)}{\log{X}} \left\{d_{3, 1}(q) + c_{3, 1} +O\left( \frac{1}{\log^2{X}} \right) \right\} + O\left(X^{1/3 + \varepsilon} \right),
\end{equation}
where
\begin{equation}
\label{eq:d3}
d_{3, 1}(q) = 6 \sum_{p | q} \frac{ \log{p}}{(p^2 - 1)(p + 1)}
\end{equation}
and
\begin{equation}
\label{eq:c3}
c_{3, 1} = -6 \sum_{p \neq 2} \frac{ \log{p}}{(p^2 - 1)(p + 1)}.
\end{equation}
The implied constants depend only on $w$, $\phi$,  and $\varepsilon$.
\end{mylemma}
Proof.  Let
\begin{equation}
U(p,\nu) = p^{-\nu/2} \mathop{ \mathop{ \mathop{ \mathop{\sum \sum}_{a \equiv r_1 \shortmod{2q}} }_{b \equiv t_1 \shortmod{2q}} }_{N_{a, b} \equiv 0 \shortmod{p} } }_{ (a, b) = 1 } \lambda_{a, b}^{\nu}(p) \, w\left( \frac{a}{A}, \frac{b}{B} \right),
\end{equation}
so that 
\begin{equation*}
\sum_{E_{a, b} \in \mathcal{F}_1} T_1(a, b) w_X\left(E_{a, b} \right)  = - \frac{2}{\log{X}} \sum_{(p, 2q) = 1} \log{p} \sum_{\nu = 1}^{\infty}  U(p,\nu) \widehat{\phi}\left( \nu \frac{\log{p}}{\log{X}}\right).
\end{equation*}
Write $U(p, \nu) = U_1 + U_2 +U_3$, say, where $U_1, U_2$, and $U_3$ correspond to $a \equiv 0 \pmod{p}$, $b \equiv 0 \pmod{p}$, and $a + 2b \equiv 0 \pmod{p}$, respectively (these conditions are mutually exclusive).  When $\nu$ is even, $\lambda_{a, b}^{\nu}(p) = p^{-\nu/2}$ for all $a$, $b$, and $p$ under consideration.  Thus for $\nu$ even and $p \neq 2$ we have (recall $\lambda(2) = 0$)
\begin{align*}
U_1 & = p^{-\nu} \mathop{ \mathop{ \mathop{ \mathop{ \mathop{ {\sum \sum} }_{a \equiv 0 \shortmod{p}}}_{a \equiv r_1 \shortmod{2q}} }_{b \equiv t_1 \shortmod{2q}}} }_{ (a, b) = 1 }  \, w\left( \frac{a}{A}, \frac{b}{B} \right) \\
& = \delta(p, 2q) p^{-\nu} \frac{1}{p + 1} M(\mathcal{F}_1) + O(p^{-\nu} X^{1/3 + \varepsilon}),
\end{align*}
where $\delta(p, 2q)$ is the indicator function of $(p, 2q) = 1$, by Lemma \ref{lem:sizelemma}.
Of course $U_2$ gives the same contribution as $U_1$.  The change of variable $a \rightarrow a - 2b$ shows $U_3$ also gives the same contribution as $U_1$ (the change of variables does not alter $\widehat{w}(0, 0)$).   By Lemma \ref{lem:lambdaformula} we have for $\nu$ odd
\begin{equation*}
U_1 = p^{-\nu} \mathop{ \mathop{ \mathop{ \mathop{ \mathop{ {\sum \sum} }_{a \equiv 0 \shortmod{p}}}_{a \equiv r_1 \shortmod{2q}} }_{b \equiv t_1 \shortmod{2q}}} }_{ (a, b) = 1 } \left(\frac{2b}{p} \right) w\left(\frac{a}{A}, \frac{b}{B} \right).
\end{equation*}
Simply apply Lemma \ref{lem:polya} to obtain the bound for $\nu$ odd
\begin{equation*}
U_1 \ll p^{-\nu - 1/2} X^{1/3 + \varepsilon}.
\end{equation*}
Again, $U_2$ and $U_3$ satisfy the same bound as $U_1$.

Gathering these estimates, we obtain
\begin{align*}
\sum_{E_{a, b} \in \mathcal{F}_1} T_1(a, b) w_X(E_{a, b}) 
& = - 6 \frac{M(\mathcal{F}_1)}{\log{X}}  \sum_{(p, 2q) = 1} \log{p} \sum_{\nu = 1}^{\infty}  \frac{\; \; p^{-2\nu} }{p + 1} \widehat{\phi}\left( 2\nu \frac{\log{p}}{\log{X}}\right) + O(X^{1/3 + \varepsilon}).
\end{align*}
Now use the approximation $\widehat{\phi}(t) = \widehat{\phi}(0) + O(t^2)$ and execute the summation over $\nu$ to get
\begin{equation*}
 - \frac{6\widehat{\phi}(0)}{\log{X}} M(\mathcal{F}_1) \left\{ \sum_{(p, 2q) = 1} \frac{\log{p}}{(p^2 - 1)(p + 1)} + O\left( \frac{1}{\log^2{X}} \right)\right\} + O(X^{1/3 + \varepsilon}),
\end{equation*}
and the proof is complete. \eop

We mention that we could have used a higher-order Taylor polynomial to approximate $\widehat{\phi}$ to obtain a better asymptotic in (\ref{eq:Tsum}).

\subsection{The Family $\mathcal{F}_2$}
\begin{mylemma} Let $T_2$ be defined as in Lemma \ref{lem:Tlemma}, but for the family $\mathcal{F}_2$.  Then
\begin{equation}
\label{eq:T2sum}
\sum_{E_{a, b} \in \mathcal{F}_1} T_2(a, b) w_X(E_{a, b}) = M(\mathcal{F}_2) \frac{\widehat{\phi}(0)}{\log{X}} \left\{d_{3, 2}(q) + c_{3, 2} + O\left(\frac{1}{ \log^2{X} } \right) \right\} + O\left(X^{1/2 + \varepsilon} \right),
\end{equation}
where
\begin{equation}
\label{eq:d3'}
d_{3, 2}(q) = 4 \sum_{p |q} \frac{ \log{p}}{(p^2 - 1)(p + 1)}
\end{equation}
and
\begin{equation}
\label{eq:c3'}
c_{3, 2} = -4 \sum_{p \neq 2} \frac{ \log{p}}{(p^2 - 1)(p + 1)}.
\end{equation}
The implied constants depend only on $w$, $\phi$,  and $\varepsilon$.
\end{mylemma}
Proof.  The proof is essentially identical to that of Lemma \ref{lem:Tlemma}.  The reason the constant differs is that the conductor of a curve in $\mathcal{F}_2$ splits into two essentially independent factors, whereas for $\mathcal{F}_1$ there are three essentially independent factors. \eop

In the course of discussion we proved
\begin{mycoro}
\label{coro:pdividingN}
Let $\delta_1 = 3$, $\delta_2 = 2$.  Then for $(p, 2q) = 1$ we have
\begin{equation*}
\mathop{ \mathop{ \sum \sum }_{E_{a, b} \in \mathcal{F}_i} }_{N_{a, b} \equiv 0 \shortmod{p}} \lambda_{a, b}^{\nu}(p) =
\begin{cases}
\delta_i M(\mathcal{F}_i) \frac{p^{-\nu/2}}{p + 1} + O(p^{-\nu/2} X^{\frac{1}{\delta_i} + \varepsilon} ) & \text{for $\nu$ even} \\
O\left( p^{-(\nu + 1)/2} X^{\frac{1}{\delta_i} + \varepsilon} \right) & \text{for $\nu$ odd}.
\end{cases}
\end{equation*}
The implied constants depend only on $w$, $\phi$, and $\varepsilon$.
\end{mycoro}

\section{Evaluating a Main Term}
\subsection{The Prime Number Theorem}
It will be necessary to use the following
\begin{mylemma} 
\label{lem:PNT}
Set $\theta(t) = \sum_{p \leq t} \log{p}$ and let 
\begin{equation}
\label{eq:R(t)}
R(t) = \theta(t) -t.
\end{equation}
Then
\begin{equation*}
\sum_{p} \frac{2\log{p}}{p \log{X}} \widehat{\phi}\left(  \frac{2 \log{p}}{\log{X}} \right) 
=
\frac{\phi(0)}{2} + \frac{2 \widehat{\phi}(0)}{\log{X}}\left(1 + \int_1^{\infty} \frac{R(t)}{t^2} dt \right) + O\left( \frac{1}{\log^3{X}} \right).
\end{equation*}
\end{mylemma}
Proof.  Let $S$ be the sum on the right hand side above.  Then by partial summation and by replacing $\theta(t)$ by $t + R(t)$ we obtain
\begin{align*}
S & = - \frac{2}{\log{X}} \int_1^{\infty} \theta(t) \left\{ - \frac{1}{t^2} \widehat{\phi}\left(  \frac{2 \log{t}}{\log{X}} \right) + \frac{2}{t^2 \log{X}} \widehat{\phi}{\,}'\left(  \frac{2 \log{t}}{\log{X}} \right) \right\} dt \\
& = \frac{2}{\log{X}} \int_1^{\infty} \frac{1}{t} \left\{ \widehat{\phi}\left(  \frac{2 \log{t}}{\log{X}} \right) 
- \frac{2}{\log{X}} \widehat{\phi}{\,}'\left(  \frac{2 \log{t}}{\log{X}} \right) \right\} dt \\
& \qquad + \frac{2}{\log{X}} \int_1^{\infty} \frac{R(t)}{t^2} \left\{ \widehat{\phi}\left(  \frac{2 \log{t}}{\log{X}} \right) 
- \frac{2}{\log{X}} \widehat{\phi}{\,}'\left(  \frac{2 \log{t}}{\log{X}} \right) \right\} dt.
\end{align*}
Applying the obvious change of variables to calculate the first integral and using the approximations $\widehat{\phi}(x) = \widehat{\phi}(0) + O(x^2)$ and $\widehat{\phi}{\,}'(x) = O(x)$ to estimate the second integral, we obtain
\begin{align*}
S & = \int_0^{\infty} \widehat{\phi}(u) du - \frac{2}{\log{X}} \int_0^\infty \widehat{\phi}\,'(u) du + \frac{2 \widehat{\phi}(0)}{\log{X}} \int_1^\infty \frac{R(t)}{t^2} dt + O\left( \frac{1}{\log^3{X}} \right) \\
& = \frac{\phi(0)}{2} + \frac{2 \widehat{\phi}(0)}{\log{X}}\left(1 + \int_1^{\infty} \frac{R(t)}{t^2} dt \right) + O\left( \frac{1}{\log^3{X}} \right),
\end{align*}
as desired. \eop
\subsection{Both Families}
Now we can easily prove
\begin{mylemma}  For $i = 1, 2$ we have
\begin{multline*}
\sum_{E_{a, b} \in \mathcal{F}_i} \sum_{p \notdiv N_{a, b}} \frac{2\log{p}}{p \log{X}} \widehat{\phi}\left(  \frac{2 \log{p}}{\log{X}} \right) w_X(E_{a, b}) \\
= M(\mathcal{F}_i)\left\{ \frac{\phi(0)}{2} + d_{4, i}(q) \frac{\widehat{\phi}(0)}{\log{X}}+ c_{4, i} \frac{\widehat{\phi}(0)}{\log{X}} + O\left(\frac{1}{\log^3{X}} \right) \right\} + O(X^{1/\delta_i + \varepsilon}),
\end{multline*}
where $d_{4, i}$ is given by \eqref{eq:d4}, $c_{4, i}$ is given by \eqref{eq:c4}, $\delta_1 = 3$ and $\delta_2 = 2$.  The implied constants depend only on $w$, $\phi$, and $\varepsilon$.
\end{mylemma}
Proof.  We do the calculation for both families simultaneously.  We first remove the restriction $p \notdivtext N_{a, b}$.  We clearly have
\begin{align*}
\sum_{E_{a, b} \in \mathcal{F}_i} \sum_{p \notdiv N_{a, b}} \frac{2\log{p}}{p \log{X}} \widehat{\phi}\left(  \frac{2 \log{p}}{\log{X}} \right) w_X(E_{a, b}) 
& = 
\sum_{E_{a, b} \in \mathcal{F}_i} w_X(E_{a, b}) \sum_{p} \frac{2\log{p}}{p \log{X}} \widehat{\phi}\left(  \frac{2 \log{p}}{\log{X}} \right) \\
& \quad - \sum_{p} \frac{2\log{p}}{p \log{X}} \widehat{\phi}\left(  \frac{2 \log{p}}{\log{X}} \right) \mathop{\sum_{E_{a, b} \in \mathcal{F}_i}}_{N_{a, b} \equiv 0 \shortmod{p}}  w \left( \frac{a}{A}, \frac{b}{B} \right).
\end{align*}
By Lemma \ref{lem:PNT} the first sum is
\begin{equation*}
M(\mathcal{F}_i) \left\{ \frac{\phi(0)}{2} + \frac{2 \widehat{\phi}(0)}{\log{X}}\left(1 + \int_1^{\infty} \frac{R(t)}{t^2} dt \right) + O\left( \frac{1}{\log^3{X}} \right) \right\}.
\end{equation*}
By Corollary \ref{coro:pdividingN} (with $\nu = 0$) we have for $(p, 2q) = 1$
\begin{equation*}
\mathop{\sum_{E_{a, b} \in \mathcal{F}_i}}_{N_{a, b} \equiv 0 \shortmod{p}}  w \left( \frac{a}{A}, \frac{b}{B} \right) = \delta_i M(\mathcal{F}_i) \frac{1}{p + 1} + O(X^{1/\delta_i + \varepsilon}),
\end{equation*}
where $\delta_1 = 3$ and $\delta_2 = 2$.  Summing over $(p, 2q) = 1$ gives
\begin{eqnarray*}
\frac{M(\mathcal{F}_i)}{\log{X}} \sum_{(p, 2q) = 1} \frac{2 \delta_i \log{p}}{p(p+1)}  \widehat{\phi}\left(  \frac{2 \log{p}}{\log{X}} \right) 
 + O(X^{1/\delta_i + \varepsilon})
\end{eqnarray*}
\begin{equation*}
= M(\mathcal{F}_i)\left\{ \frac{2 \delta_i \widehat{\phi}(0)}{\log{X}} \sum_{(p,2q) = 1}  \frac{\log{p}}{p(p+1)} + O\left(\frac{1}{\log^3{X}} \right)   \right\} + O(X^{1/\delta_i + \varepsilon}).
\end{equation*}
The contribution from $p = 2$ is simply
\begin{equation*}
M(\mathcal{F}_i) \left\{ \frac{\log{2} }{\log{X}}\widehat{\phi}(0)  + O\left(\frac{1}{\log^3{X}} \right)   \right\}.
\end{equation*}

By gathering terms the proof of the lemma is complete with
\begin{align}
\label{eq:d4}
d_{4, i} &= 2 \delta_i \sum_{p | q} \frac{\log{p}}{p(p+1)}, \\
\label{eq:c4}
c_{4, i} &=  2 \left(1 + \int_1^{\infty} \frac{R(t)}{t^2} dt  - \delta_i \sum_{p \neq 2} \frac{\log{p}}{p(p+1)} - \log{2} \right),
\end{align}
with $R(t)$ given by (\ref{eq:R(t)}) (also recall $\delta_1 = 3$, $\delta_2 = 2$).
\eop

\section{The Variation From $a_0$ and $b_0$}
\label{section:variation}
\subsection{Notation}
\begin{mydefi}
Set
\begin{equation*}
\varphi(p, \nu) =  \widehat{\phi}\left(  \frac{\nu \log{p}}{\log{X}} \right) - p^{-1} \widehat{\phi}\left(  \frac{(\nu + 2) \log{p}}{\log{X}} \right).
\end{equation*}
For $\mathcal{F}_i$, $i = 1, 2$, set
\begin{equation*}
Q_i(p^\nu) = \mathop{ \mathop{\sum \sum}_{\alpha \shortmod{p}}}_{\beta \shortmod{p}} \lambda_{\alpha, \beta}(p^{\nu}).
\end{equation*}
Further, set
\begin{equation}
\label{eq:fourierpart}
J(\mathcal{F}_i) = 
- 2 \sum_{E_{a, b} \in \mathcal{F}_i}  \sum_{p \notdiv N_{a, b}}  \frac{  \log{p}}{\log{X}}  \sum_{\nu = 1}^{\infty} \frac{\lambda_{a, b}(p^\nu)}{p^{\nu/2}} 
\varphi(p, \nu)
w_X(E_{a, b}).
\end{equation}
\end{mydefi}
For notational cleanliness we do not exhibit the dependence of $N_{a, b}$ and $\lambda_{a, b}(p^\nu)$ on the family $\mathcal{F}_i$.
\subsection{A Complete Character Sum Computation}
\label{section:complete}
The character sums $Q_i(p^\nu)$ hold important arithmetical information about the distribution of zeros in families of elliptic curves.  In this section we show $Q_i(p^{\nu}) = 0$ for $\nu = 2$ and all $\nu$ odd.  It is of interest to study such character sums for other families of elliptic curves.  When obtaining a density theorem with large support for a family of elliptic curves it is necessary to study variants of $Q_i(p)$ where the sum is twisted by an arbitrary additive character modulo $p$.  Such investigations are undertaken in \cite{Young}.

We first show that $Q_1(p^{\nu}) = 0$ for $\nu$ odd.  It is easy to see that for $(p, N) = 1$ we have $\lambda(p^{\nu}) = \sum_l d_l \lambda^l(p)$, where $d_l = 0$ if $l \not \equiv \nu \pmod{2}$.  Hence we have
\begin{equation*}
\sum_{\alpha} \sum_{\beta} \lambda_{\alpha, \beta}(p^\nu) = \sum_{l \equiv \nu \shortmod{2}} d_l'\sum_\alpha \sum_\beta \lambda_{\alpha, \beta}^l(p) + \sum_{l \equiv \nu \shortmod{2}} d_l'' \mathop{\sum_\alpha \sum_\beta}_{p | N_{\alpha, \beta}} \lambda_{\alpha, \beta}^l(p),
\end{equation*}
for some $d_l', d_l'' \in \mc$.
The problem reduces to showing each of the two sums above vanish.

For the first sum, we see that
\begin{equation*}
p^{l/2} \sum_{\alpha} \sum_{\beta} \lambda_{\alpha, \beta}^{l}(p) = \mathop{\sum \cdots \sum}_{y_1, \ldots, y_l} \sum_{\alpha} \sum_{\beta} \prod_{i=1}^{l} \left(\frac{y_i(y_i - \alpha)(y_i + 2 \beta)}{p} \right).
\end{equation*}
Take $e$ such that $(e/p) = -1$ and apply the change of variables $y_i \rightarrow e y_i$, $\alpha \rightarrow e \alpha$, $\beta \rightarrow e \beta$.  We get the same sum multiplied by $(e/p)^{3l}$; hence the sum is zero for $l$ odd.  A glance at (\ref{eq:lambdaformular}) shows that the second sum vanishes.

The argument showing $Q_2(p^\nu) = 0$ for $\nu$ odd is similar.

We now show $Q_1(p^2)=0$.  We use the identity $\lambda(p^2) = \lambda^2(p) - 1$ for $(p, N) = 1$ and $\lambda(p^2) = \lambda^2(p)$ for $p | N$.  The derivation is simple (assume $p \neq 2$):
\begin{eqnarray*}
Q_1(p^2) & = & \sum_\alpha \sum_\beta \lambda_{\alpha, \beta}^2(p)  - \mathop{\sum_{\alpha} \sum_{\beta}}_{\alpha \beta (\alpha + 2\beta) \not \equiv 0 \shortmod{p}} 1 \\
& = & \frac{1}{p} \sum_x \sum_y \sum_\alpha \sum_\beta  \left(\frac{x(x - \alpha)(x + 2 \beta)}{p} \right) \left(\frac{y(y - \alpha)(y + 2 \beta)}{p} \right) - (p-1)(p-2) \\
& = & \frac{1}{p} \left[\mathop{\sum \sum}_{x \neq y} \left(\frac{xy}{p}\right) + (p-1)^2\sum_x \left(\frac{x^2}{p}\right) \right] - (p^2 - 3p + 2) \\
& = &  \frac{1}{p} \left[-p + 1  + (p-1)^3\right] - (p^2 - 3p + 2) \\
& = & 0.
\end{eqnarray*}

Similarly,
\begin{eqnarray*}
Q_2(p^2) & = & \sum_\alpha \sum_\beta \lambda_{\alpha, \beta}^2(p) - \mathop{\sum_{\alpha} \sum_{\beta}}_{\beta(\alpha^2 + \beta) \not \equiv 0 \shortmod{p}} 1  \\
& = &   \frac{1}{p} \sum_x \sum_y \sum_\alpha \sum_\beta  \left(\frac{x(x^2 + 2\alpha x - \beta)}{p} \right) \left(\frac{y(y^2 + 2 \alpha y - \beta)}{p} \right) -(p -1)^2\\
& = &  \frac{1}{p} \sum_x \sum_y \sum_\alpha \sum_\beta  \left(\frac{xy}{p}\right) \left(\frac{\beta}{p}\right) \left(\frac{\beta + (y^2 - x^2) + \alpha(y - x)}{p}\right) - (p-1)^2\\
& = &   \sum_x \sum_\beta \left(\frac{x^2}{p}\right) \left(\frac{\beta^2}{p}\right) -(p-1)^2\\
& = & 0.
\end{eqnarray*}

Note that $Q_i(p^k)$ is not always zero.  For instance, $Q_1(5^4) = -216/25$, $Q_1(7^4) = 528/49$, and $Q_2(3^6) = -8/9$.  One can show that $Q_2(p^4) = 0$ for all $p$, but the calculations are lengthy and tangential to our purpose here.

\subsection{Both Families}
\begin{mylemma}
\label{lem:J}
We have
\begin{multline*}
J(\mathcal{F}_i) = M(\mathcal{F}_i) \frac{\widehat{\phi}(0)}{\log{X}} \left\{c_{5, i} + c_{6, i}  + d_{5, i} (q) + d_{6,i}(q) + e_i(a_0, b_0) +O\left(\frac{(\log{q})^{1/2 + \varepsilon}}{\log^2{X}}\right) \right\} \\
+ O((AB)^{1-\varepsilon}),
\end{multline*}
where
\begin{eqnarray}
\label{eq:c5}
c_{5, i} & = & -2   \sum_{p \neq 2} \frac{\log{p}}{p + 1} \sum_{l = 1}^{\infty} \frac{Q_i(p^{2l})}{p^{l + 1}}, \\
\label{eq:c6}
c_{6, i} & = & 2 \delta_i \sum_{p \neq 2} \frac{\log{p}}{p(p + 1)^2} \\
\label{eq:d5}
d_{5, i} (q)  & = &     2 \sum_{p | q} \frac{\log{p}}{p + 1} \sum_{l = 1}^{\infty} \frac{Q_i(p^{2l})}{p^{l + 1}},  \\
\label{eq:d6}
d_{6, i}(q) & = & - 2 \delta_i \sum_{p | q} \frac{\log{p}}{p(p+1)^2} \\
\nonumber \text{and} \\
\label{eq:e1}
e_i(a_0, b_0) & = & -  2 \sum_{p | q} \log{p} \left( 1- \frac{1}{p} \right)  \sum_{\nu = 1}^{\infty} \frac{\lambda_{a_0, b_0}(p^{\nu})}{p^{\nu/2}}.
\end{eqnarray}
The implied constants depend only on $w$, $\phi$, and $\varepsilon$.  Further, $\varepsilon$ must be sufficiently small with respect to the support of $\widehat{\phi}$.
\end{mylemma}
Recall $\delta_1 = 3$ and $\delta_2 = 2$.  We can also use the identity
\begin{equation*}
\sum_{\nu = 0}^{\infty} \frac{\lambda(p^{\nu})}{p^{\nu/2}} = \left(1 - \frac{\lambda(p)}{p^{1/2}} + \frac{\chi_N(p)}{p} \right)^{-1}
\end{equation*}
to find alternate expressions for $c_{5, i}$, $d_{5, i}(q)$, and $e_i(a_0, b_0)$.  Here $\chi_N$ is the principal Dirichlet character $\mymod{N}$.  Namely we have
\begin{align}
\tag{\ref{eq:c5}'}
c_{5, i} & =  -2  \sum_{p \neq 2} \frac{\log{p}}{p(p + 1)} \mathop{\sum \sum}_{a, b \shortmod{p}} \left\{ \left(1 - \frac{\lambda_{a, b}(p)}{p^{1/2}} + \frac{\chi_{N_{a,b}}(p)}{p} \right)^{-1} - 1\right\}, \\
\tag{\ref{eq:d5}'}
d_{5, i} (q)  & =   2 \sum_{p | q} \frac{\log{p}}{p(p + 1)} \mathop{\sum \sum}_{a, b \shortmod{p}} \left\{ \left(1 - \frac{\lambda_{a, b}(p)}{p^{1/2}} + \frac{\chi_{N_{a,b}}(p)}{p} \right)^{-1} - 1\right\},     \\
\nonumber \text{and} \\
\tag{\ref{eq:e1}'}
e_i(a_0, b_0) & =  -  2  \sum_{p | q} \log{p} \left( 1- \frac{1}{p} \right) \left\{ \left(1 - \frac{\lambda_{a_0, b_0}(p)}{p^{1/2}} + \frac{1}{p} \right)^{-1} - 1 \right\}.
\end{align}
It is interesting that $L_p(1/2, E_{a_0, b_0})$ is visibly present in the expression for $e_i$.

Proof.  The basic idea is to simply complete the sum over $a$ and $b$ modulo $p$ to write our sum in terms of a sum of $Q(p^{\nu})$.  For $\nu \geq 3$ this is completely straightforward.  For $\nu =2$ some work needs to be done to show that the summation over $p$ converges.  The calculations in the previous section showing $Q(p^2) = 0$ are for this purpose (the trivial bound $Q(p^2) \ll p^2$ barely fails to succeed in showing the sum converges).  The case $\nu = 1$ is in some sense a separate issue.  If $\widehat{\phi}$ has small support (i.e. for $p$ rather small) it is easy to show that there is no contribution from terms with $\nu = 1$.  For larger support it becomes a difficult problem to prove the same result.  See \cite{Young} for a detailed investigation into such problems for a variety of families of elliptic curves.

It is convenient to extend the summation over all primes $p$.  Write $J = J' - J''$, where $J'$ is defined by the sum (\ref{eq:fourierpart}) but with no restriction on $p$, and $J''$ is defined similarly but with the condition $p | N$.  Since $\lambda(2) =0$ we may freely assume $p \neq 2$.  

We first show
\begin{equation}
\label{eq:J''calc}
J''(\mathcal{F}_i) =  M(\mathcal{F}_i) \frac{\widehat{\phi}(0)}{\log{X}}\left( - c_{6, i} - d_{6, i} + O\left( \frac{1}{\log^2{X}} \right) \right) +O\left(X^{\frac{1}{\delta_i} + \varepsilon}\right).
\end{equation}
Using Corollary \ref{coro:pdividingN} and the approximation 
\begin{equation}
\label{eq:phiapprox}
\varphi(p, \nu) = \widehat{\phi}(0) \left(1 - \frac{1}{p}\right) + O\left(\frac{\nu^2 \log^2{p}}{ \log^{2}{X}} \right)
\end{equation}
we have
\begin{align}
\nonumber
J''(\mathcal{F}_i) & =  -2 \sum_{(p, 2q) = 1} \sum_{\nu = 1}^{\infty} \mathop{\sum_{E_{a, b} \in \mathcal{F}_i} }_{N_{a, b} \equiv 0 \shortmod{p} } \frac{\log{p}}{\log{X}} \frac{\lambda_{a, b}(p^\nu)}{p^{\nu/2}} \varphi(p, \nu) w \left( \frac{a}{A}, \frac{b}{B} \right)  \\
\nonumber
& =  -2 \delta_i M(\mathcal{F}_i) \sum_{(p, 2q) = 1} \sum_{\nu \text{ even}} \frac{\log{p}}{\log{X}}  \frac{p^{-\nu}}{p + 1} \varphi(p, \nu) + O\left(X^{\frac{1}{\delta_i} + \varepsilon}\right) \\
\nonumber
& =  -2 \delta_i M(\mathcal{F}_i) \left\{ \sum_{(p, 2q) = 1} \sum_{\nu \text{ even}} \frac{\log{p}}{\log{X}}  \frac{p^{-\nu}}{p + 1} \left(1 - \frac{1}{p} \right) \widehat{\phi}(0) + O\left( \frac{1}{\log^2{X}} \right) \right\} + O\left(X^{\frac{1}{\delta_i} + \varepsilon}\right)\\
\label{eq:J''sum}
& =  -2 \delta_i M(\mathcal{F}_i) \frac{\widehat{\phi}(0)}{\log{X}} \left\{ \sum_{(p, 2q) = 1} \frac{\log{p}}{p(p + 1)^2} + O\left( \frac{1}{\log^2{X}} \right) \right\} + O\left(X^{\frac{1}{\delta_i} + \varepsilon}\right),
\end{align}
as desired.

It suffices to show 
\begin{equation*}
J'(\mathcal{F}_i) = M(\mathcal{F}_i) \frac{\widehat{\phi}(0)}{\log{X}} \left(c_{5, i} + d_{5, i}(q) + e_i(a_0, b_0) + O\left(\frac{1}{\log^2{X}}\right) \right).
\end{equation*}
Let
\begin{equation*}
V(p, \nu) = \mathop{ \mathop{ \mathop{\sum \sum}_{a \equiv r_i \shortmod{2^i q}} }_{b \equiv t_i \shortmod{2^i q}} }_{(a, b) = 1} \lambda_{a, b}(p^\nu) w\left( \frac{a}{A}, \frac{b}{B} \right),
\end{equation*}
so that
\begin{equation*}
J'(\mathcal{F}_i) = -2 \sum_{p} \frac{\log{p}}{\log{X}} \sum_{\nu=1}^{\infty} \frac{V(p,\nu)}{p^{\nu/2}}  \varphi(p, \nu).
\end{equation*}
We will calculate $V(p,\nu)$.  
In case $p | q$ we then have
\begin{equation*}
V(p,\nu) =  \lambda_{a_0, b_0}(p^\nu)  M(\mathcal{F}_i) + O(X^{1/\delta_i + \varepsilon}),
\end{equation*}
by Lemmas \ref{lem:mass} and \ref{lem:mass2}.
The contribution to $J'$ from $p | q$ is thus
\begin{align}
\label{eq:J'p=q}
& -2 \frac{\widehat{\phi}(0)}{\log{X}} M(\mathcal{F}_i)\sum_{p | q} \log{p}   \sum_{\nu = 1}^{\infty} \frac{ \lambda_{a_0, b_0}(p^{\nu})}{p^{\nu/2}} \varphi(p, \nu) 
+ O(X^{1/\delta_i + \varepsilon}) \\
\nonumber
& = -2 \frac{\widehat{\phi}(0)}{\log{X}} M(\mathcal{F}_i) \left\{ \sum_{p | q} \log{p} \left( 1- \frac1p \right)   \sum_{\nu = 1}^{\infty} \frac{ \lambda_{a_0, b_0}(p^{\nu})}{p^{\nu/2}} +O\left(\frac{\alpha(q)}{\log^2{X}} \right) \right\} + O(X^{1/\delta_i + \varepsilon}),
\end{align}
where 
\begin{equation*}
\alpha(q) = \sum_{p | q}\frac{\log^3{p}}{p^{1/2}} \ll (\log{q})^{1/2} \log^2(\log{3q}).
\end{equation*}
This gives the term $e_i(a_0, b_0)$.
From now on assume $p \notdivtext 2q$.  

We shall now prove the following approximation
\begin{equation}
\label{eq:Vpnutrivialestimate}
V(p,\nu) = \frac{Q_i(p^\nu)}{p^2 - 1} M(\mathcal{F}_i) + O\left(\left(\nu + 1\right)( p (A + B)^{1 + \varepsilon} + p^2 (AB)^{\varepsilon} + (AB)^{1-\varepsilon}) \right),
\end{equation}
uniform with respect to $q$.  This bound is nontrivial for $p \ll (\text{min}( A, B))^{1 - \varepsilon}$.  

Proof of (\ref{eq:Vpnutrivialestimate}).  By M\"{o}bius inversion,
\begin{equation*}
V(p,\nu) = \sum_{(d, 2pq) = 1} \mu(d)  V_d(p,\nu),
\end{equation*}
where
\begin{equation*}
V_d(p,\nu) = \mathop{ \mathop{\sum \sum}_{ad \equiv r_i \shortmod{2^i q}} }_{bd \equiv t_i \shortmod{2^i q}}  \lambda_{ad, bd}(p^\nu) w\left( \frac{ad}{A}, \frac{bd}{B} \right).
\end{equation*}
The restriction $(d, p) = 1$ can be imposed because $\lambda_{0, 0} (p^{\nu}) = 0$ for all $p, \nu$.
Since $V_d(p,\nu) \ll d^{-2} (\nu + 1) AB $ we have
for any $R \geq 1$ the trivial estimate
\begin{equation}
\label{eq:Vdtrivial}
\sum_{d \geq R} |V_d(p,\nu)| \ll \frac{ AB(\nu + 1)}{R},
\end{equation}
the implied constant depending on $w$ only.
Performing Poisson summation in $a$ and $b \pmod{2^ipq}$ gives
\begin{equation}
\label{eq:Vdformula}
V_d(p,\nu) = \frac{AB}{4^i d^2p^2 q^2}  \sum_h \sum_k Y(h, k, p, \nu) \widehat{w}\left( \frac{hA}{2^i dpq}, \frac{kB}{2^i dpq} \right),
\end{equation}
where
\begin{equation*}
Y(h, k, p, \nu) = \mathop{\mathop{\sum_{\alpha \shortmod{2^i pq}} \sum_{\beta \shortmod{2^i pq}}}_{\alpha d \equiv r_i \shortmod{2^i q}} }_{\beta d \equiv t_i \shortmod{2^i q}} \lambda_{\alpha d, \beta d}(p^\nu) e\left(\frac{\alpha h + \beta k}{2^i pq} \right).
\end{equation*}
Computing $Y$ by separating the variables via the Chinese remainder theorem we get
\begin{equation*}
Y(h, k, p, \nu) = e\left( \frac{\overline{dp}(r_i h + t_i k) }{2^i q} \right) 
\mathop{ \mathop{\sum \sum}_{\alpha \shortmod{p}}}_{\beta \shortmod{p}} \lambda_{\alpha , \beta }(p^\nu) e\left(\frac{\overline{2^i d q}(\alpha h + \beta k)}{p} \right).
\end{equation*}
If we estimate the nonzero frequencies of $V_d$ trivially we get the estimates (for $(d, 2pq) = 1$)
\begin{equation*}
V_d(p,\nu) = \frac{AB}{4^i d^2p^2q^2} \left\{ Q_i(p^\nu) + O\left(\frac{dp^{3}q (\nu+1)}{\text{min}(A, B)} + \frac{d^2p^{ 4}q^2(\nu + 1)}{AB} \right) \right\}
\end{equation*}
and hence
\begin{equation}
\label{eq:Vdsmalld}
\mathop{\sum_{d < R}}_{(d, 2pq) = 1} \mu(d) V_d(p,\nu) = \frac{Q_i(p^\nu)}{p^2 - 1} M(\mathcal{F}_i) + O\left( \left(\nu + 1\right) \left( p (A + B) \log{R} + p^2R + \frac{AB}{R} \right)\right). 
\end{equation}
Gathering the estimates (\ref{eq:Vdtrivial}) and (\ref{eq:Vdsmalld}) and taking $R = (AB)^{ \varepsilon}$ we obtain (\ref{eq:Vpnutrivialestimate}).

The case $\nu = 1$ requires special care so we first estimate the contribution to $J'(\mathcal{F}_i)$ from $\nu \geq 2$.  For convenience we denote
\begin{eqnarray*}
J_1'(\mathcal{F}_i) & = &  -2 \sum_{(p, q) = 1} \frac{\log{p}}{\log{X}} \frac{V(p,1)}{p^{1/2}}  \varphi(p, 1), \\
J_2'(\mathcal{F}_i) & = & -2 \sum_{(p, q) = 1} \frac{\log{p}}{\log{X}} \sum_{\nu=2}^{\infty} \frac{V(p,\nu)}{p^{\nu/2}}  \varphi(p, \nu).
\end{eqnarray*}
At this point take $0 < \rho_1 < 2/3$, $0 < \rho_2 < 1/2$ such that for each $\mathcal{F}_i$ we have $\text{supp } \widehat{\phi} \subset (- \rho_i, \rho_i)$.  Then $\varphi(p, \nu) = 0$ for $p^{\nu} \geq X^{\rho_i}$.  The inequality $p \ll X^{\rho_i/\nu}$ implies $p \ll (\text{min}( A, B))^{1 - \varepsilon}$ for $\nu \geq 2$, i.e. the estimate (\ref{eq:Vpnutrivialestimate}) is nontrivial in this range.  Hence we easily obtain
\begin{eqnarray}
\nonumber
J_2'(\mathcal{F}_i) & = & -2 \mathop{\sum_{(p, q) = 1}  \sum_{\nu=2}^{\infty}}_{p^{\nu} \leq X^{\rho_i}} \frac{\log{p}}{\log{X}} \frac{V(p,\nu)}{p^{\nu/2}}  \varphi(p, \nu) \\
\nonumber
& = & - \frac{2}{\log{X}} \sum_{(p, q) = 1}  \sum_{\nu=2}^{\infty} \frac{\log{p}}{p^{\nu/2}} \varphi(p, \nu) 
\frac{Q_i(p^\nu)}{p^2 - 1} M(\mathcal{F}_i) 
 + O((AB)^{1-\varepsilon}) \\
\nonumber
& = & - \frac{2 M(\mathcal{F}_i)}{\log{X}} \sum_{(p, q) = 1}  \sum_{\nu=1}^{\infty} \frac{\log{p}}{p^{2\nu}} \varphi(p, 2\nu) 
\frac{Q_i(p^\nu)}{p^2 - 1} 
 + O((AB)^{1-\varepsilon}) \\
\nonumber
& = & - \frac{2 M(\mathcal{F}_i) \widehat{\phi}(0) }{\log{X}} \left\{ \sum_{(p, q) = 1}  \sum_{\nu=1}^{\infty} \frac{\log{p}}{p^{\nu + 1}}
\frac{Q_i(p^{2\nu})}{p + 1} 
 + O\left(\frac{1}{\log^2{X}}\right) \right\} + O((AB)^{1-\varepsilon}) \\
\label{eq:J2'}
& = & M(\mathcal{F}_i) \frac{\widehat{\phi}(0) }{\log{X}} \left( c_{5, i} + d_{5, i}(q) + O\left( \frac{1}{\log^2{X}} \right) \right) + O((AB)^{1-\varepsilon}).
\end{eqnarray}
Here $\varepsilon$ is bounded from above in terms of the support of $\widehat{\phi}$.  As $\rho_1$ approaches $2/3$ or $\rho_2$ approaches $1/2$ we must have $\varepsilon$ approaching zero.

It remains to show
\begin{equation*}
J_1'(\mathcal{F}_i) \ll (AB)^{1- \varepsilon}
\end{equation*}
for test functions $\widehat{\phi}$ supported in the desired ranges.  The proofs (distinct for each family) are essentially carried out in \cite{Young}.  For small support the result follows from (\ref{eq:Vpnutrivialestimate}).

By gathering the terms (\ref{eq:J''sum}), (\ref{eq:J'p=q}), and (\ref{eq:J2'}), the proof of Lemma \ref{lem:J} is now complete. \eop

\end{document}